\documentclass{article}

\usepackage{PRIMEarxiv}

\usepackage[utf8]{inputenc} % allow utf-8 input
\usepackage[T1]{fontenc}    % use 8-bit T1 fonts
\usepackage[colorlinks,bookmarksopen,bookmarksnumbered,citecolor=blue,urlcolor=blue]{hyperref}
\usepackage{url}            % simple URL typesetting
\usepackage{booktabs}       % professional-quality tables
\usepackage{amsfonts}       % blackboard math symbols
\usepackage{nicefrac}       % compact symbols for 1/2, etc.
\usepackage{microtype}      % microtypography
\usepackage{lipsum}
\usepackage{fancyhdr}       % header
\usepackage{graphicx}       % graphics
\graphicspath{{media/}}     % organize your images and other figures under media/ folder

\usepackage{amssymb}
\usepackage[tbtags]{amsmath}
\usepackage{amsfonts}
\usepackage{amsbsy}
\usepackage{mathrsfs}
\usepackage{esint}

\usepackage{bm}
\usepackage{subeqnarray}
\usepackage{cases}
\usepackage{subfigure}
\usepackage{booktabs, multirow}
\usepackage{enumerate}
\usepackage{graphicx}
\usepackage{threeparttable}
\usepackage[utf8]{inputenc}
\usepackage[T1]{fontenc}
\usepackage{lmodern}
\usepackage[figurename=Fig.,labelfont=bf,labelsep=period]{caption}
\usepackage{appendix}

\title{A multi-scale framework for neural network enhanced methods to the solution of partial differential equations
}

\author{
  Xiaodan Ren\\
  College of Civil Engineering\\
  Tongji University\\
  Shanghai, China\\
  \texttt{rxdtj@tongji.edu.cn} \\
}

\begin{document}
\maketitle

\begin{abstract}
In the present work, a multi-scale framework for neural network enhanced methods is proposed for approximation of function and solution of partial differential equations (PDEs). By introducing the multi-scale concept, the total solution of the target problem could be decomposed into two parts, i.e. the coarse scale solution and the fine scale solution. In the coarse scale, the conventional numerical methods (e.g. finite element methods) are applied and the coarse scale solution could be obtained. In the fine scale, the neural networks is introduced to formulate the solution. The custom loss functions are developed by taking into account the governing equations and boundary conditions of PDEs, the constraints and the interaction from coarse scale. The proposed methods are illustrated and examined by various of testing cases.
\end{abstract}

% keywords can be removed

%% main text

\section{Introduction}

In the last 10-20 years, artificial neural networks have been experiencing explosive advances\cite{Goodfellow-et-al-2016} and enforced the rapid developments in numbers of scientific and engineering disciplines including computer vision\cite{Lecun1995, Krizhevsky2012}, natural language processing\cite{Cho2014,Wu2016}, bio-technology\cite{Alipanahi2015} and so on. In the field of scientific computing, numbers of developments have been made in solving differential equations with the help of artificial neural networks \cite{Yadav2015,Long2018,LONG2019108925,khoo_lu_ying_2021}. Among them, the methods that trained the neural networks based on unsupervised procedures with the custom loss functions derived from the differential equations and boundary conditions gained enormous attention in the last few years. E and Yu \cite{EWeinan2018} proposed the Deep Ritz Method, which used the neural network representation of functions in the context of the Ritz Method. In their work, the validation and effectiveness of the Deep Ritz Method were shown via the numerical solutions of Poisson equations in low and high dimensions. Later, Liao and Ming \cite{Ming_2021} extended the Deep Ritz Method to Deep Nitsche Method by imposing the essential boudary conditions with the Nitsche's method. The physics-informed neural network (PINN) was proposed by Raissi and coworkers \cite{Raissi2018932,RAISSI2019686}. PINN targeted the solution of PDEs in the strong form and lead to a numerical procedure in the form collocation methods. It was shown by a series of numerical examples that the discontinuity for the solution of Burger's equation and the vortices for the solution of Navier–Stokes equation could be captured by PINN. Then the  Galerkin weak form was introduced to PINN and the variational physics-informed neural network  (VPINN) was developed \cite{EhsanKharazmi2019, VarNet2019}. Furthermore, the sub-domian Petrov-Galerkin is introduced to PINN and the \textit{hp}-Variational Physics Informed Neural Network (\textit{hp}-VPINN) is developed\cite{KHARAZMI2021113547}. This method offers \textit{hp}-refinements based on domain decomposition and high-order polynomial test functions. Samaniego \textit{et. al.}\cite{SAMANIEGO2020112790} also demonstrated that the energy based loss functions worked well with less number of unknowns than the collocation based methods. Shin \textit{et. al.} \cite{Shin2020} discussed the convergence of PINN for linear second order PDEs. Karniadakis \textit{et.~al.} \cite{Karniadakis2021} reviewed the recent progresses related to physics informed machine learning methods.

A PINN architecture based Python package named SciANN for scientific computing was developed by Haghighat and Juanes \cite{HAGHIGHAT2021113552}. Based on SciANN, the problems of solid structure undergoing plasticity was investigated and the strain localization could be captured. On the other hand, the total number of parameters for the neural network used for the problem was about 100 millions. That is to say, the NN with densely connected hidden layers was over-parameterized and of relatively low efficiency. Zhang and coworkers \cite{Zhang2021} designed the structure of NN based on conventional numerical methods such as finite element methods and reproducing kernel (RK) particles. Due to the sparse nature of these methods, the designed networks were sparse and the numbers of connections and parameters may be reduced. Recently, Baek \textit{et. al.} \cite{Baek2022} focused on the strain localization problem of solid and proposed a neural network-enhanced reproducing kernel particle method.

As can be concluded from the existing works, the neural networks approximate the solutions of differential equations in a nonlinear construction\cite{devore_1998,devore_2021}. The composition of hidden layers offers highly flexible function spaces for solution searching. 
On the other hand, the conventional studies of scientific and engineering computing often established the approximation of functions in the way of additive construction\cite{EWeinan2018}, which was linear in its nature. Then all kinds of classic methods, e.g. finite element methods\cite{Hughes2012,Belytschko2013}, finite difference methods\cite{Ames1977,Shu1998}, meshfree methods\cite{Liu1995,CHEN1996195}, isogeometric analysis\cite{HUGHES20054135} and so on, could be developed. As can be seen, the nonlinear construction is more flexible but less efficient. Thus the NN based methods could be used to capture certain complex features that can not be easily represented by the conventional methods. On the other side, the regular solution of the PDEs could be well represented by the linear construction with reasonable efficiency and robustness. Although NN based solutions of PDEs have gained more and more attentions, the researches focused on development of methods which could make full use of NN methods and conventional methods are limited. 

In the present work, a multi-scale framework for neural network enhanced methods is proposed for function approximation and solution of partial differential equations. By introducing the multi-scale concept to the target problem, the solution could be decomposed into two parts, i.e. the coarse scale solution and the fine scale solution. In the coarse scale, the target problem is solved based on conventional methods, in most cases, finite element methods. In the fine scale, the solution is developed base on neural networks. The governing equation, boundary conditions, constraints and the interaction from the coarse scale solutions are implemented in the loss functions. The article is organized as follows. In Section \ref{sec:II}, the formulation of NN enhanced approximation will be given as the basis for numerical solution of differential equations. In Section \ref{sec:III}, the general multi-scale representation of variational form of PDEs is given and the solution strategies are discussed. The neural network enhanced methods and the numerical algorithm are developed in Section \ref{sec:IV}. In Section \ref{sec:V}, numbers of typical cases are given. This is followed by the concluding remarks in Section \ref{sec:VI}.

\section{Neural network enhanced approximation}
\label{sec:II}

\subsection{Multi-scale NN enhanced approximation}

For a given function $y = f(x)$, the conventional methods of approximation could be expressed in the additive form as follows
\begin{equation}\label{eq:kernel_int}
y(\bm{x}) \approx y^\text{h}(\bm{x}) = \sum_{J = 1}^{N_\text{P}} \Psi_J(\bm{x}) d_J
\end{equation} 
The approximate/shape/kernel functions $\Psi_J(\bm{x})$ could be established based on some conventional methods. If the functions $\Psi_J(\bm{x})$ satisfy the property of Kronecker delta, the coefficients should be the value of function, e.g. $d_J = y_J = y(\bm{x}_J)$.

We consider the multi-scale form of approximation as follows
\begin{equation}
y(\bm{x}) \approx y^\text{h}(\bm{x}) = \bar{y}^\text{h}(\bm{x}) + \tilde{y}^\text{h}(\bm{x})
\end{equation}
where $\bar{y}^\text{h}(\bm{x})$ and $\tilde{y}^\text{h}(\bm{x})$ are approximation functions in coarse-scale and fine-scale, respectively. Recall the conventional methods shown in Eq.(\ref{eq:kernel_int}), the coarse-scale approximation function is expressed in the following form
\begin{equation}
\bar{y}^\text{h}(\bm{x}) = \sum_{J = 1}^{N_\text{P}} \Psi_J(\bm{x}) d_J
\end{equation}
In the meanwhile, it is proposed that the fine-scale approximation function could be established in the form of neural networks. Usually, NN is defined based on totally $l$ hidden layers. Thus we have
\begin{equation}\label{eq:NN_module}
\tilde{y}^\text{h}(\bm{x}) = \mathcal{N}(\bm{x};\bm{W},\bm{b}) =
\mathcal{G} \circ \Phi_l \circ \Phi _{l-1} \circ \cdots \circ \Phi_1 (\bm{x})
\end{equation}
For the $i$-th hidden layer, it has $N^{\text{HL}}_i$ neurons with the activation function in the following form
\begin{equation}\label{eq:NN_layer}
\Phi_i (\bm{z}) = \sigma_i (\bm{W}_i \times \bm{z} + \bm{b}_i) 
\end{equation}
where $\bm{W}_i \in \mathbb{R}^{N^{\text{HL}}_{i-1} \times N^{\text{HL}}_i}$ and $\bm{b}_i \in \mathbb{R}^{N^{\text{HL}}_i}$ are the weight matrix and the bias vector for the $i$-th hidden layer, respectively. For the input layer (the zeroth layer), $N^{\text{HL}}_0$ is the dimension of input defined by the problem. For the output layer $\mathcal{G}(\cdot)$, the dimension $N^{\text{HL}}_{i+1}$ is the output dimension defined by the problem and the linear mapping is usually adopted. Thus we have
\begin{equation}
\mathcal{G} (\bm{z}) = \bm{W}_{l+1} \times \bm{z} + \bm{b}_{l+1}
\end{equation}
The parameters in the neural network $\mathcal{N}(\bm{x};\bm{W},\bm{b})$ are $\bm{W} = [\bm{W}_1, \bm{W}_2,...,\bm{W}_l,\bm{W}_{l+1}]$ and $\bm{b} = [\bm{b}_1,\bm{b}_2,...,\bm{b}_l,\bm{b}_{l+1}]$. It can be seen that the relation between the target function $\tilde{y}^\text{h}(\bm{x})$ and the parameters of the approximation $\bm{W}$ and $\bm{b}$ is highly nonlinear. Thus we quote from \cite{devore_1998,devore_2021} that NN is a kind of nonlinear approximate method.

How to obtain the parameters based on proper training process is critical for this problem. First consider the loss function in the general form. We have
\begin{equation}\label{eq:multiscale_loss_function}
\begin{split}
\mathbf{Loss}
&= \|y(\bm{x}) - y^\text{h}(\bm{x})\|^2_{L^2} \\
&=\|y(\bm{x}) - \bar{y}^\text{h}(\bm{x}) - \tilde{y}^\text{h}(\bm{x})\|^2_{L^2} \\
&=\left\|y(\bm{x}) -\sum_{J = 1}^{N_\text{P}} \Psi_J(\bm{x}) d_J - \mathcal{N}(\bm{x};\bm{W},\bm{b}) \right\|^2_{L^2}
\end{split}
\end{equation}
where the $L^2(\Omega)$ norm of function $f(\bm{x})$ defined on a bounded domain $\Omega \subset \mathbb{R}^d$ is expressed as follows
\begin{equation}
\|f(\bm{x}) \|^2_{L^2} = (f(\bm{x}),f(\bm{x})) = \int_{\Omega} f^2(\bm{x}) \text{d} \Omega
\end{equation}
With minimization of loss function, the values of approximate coefficients $\bm{d} = (d_1, d_2, ....,d_{N_\text{P}})^{\text{T}}$ and the unknowns $(\bm{W},\bm{b})$ within the neural network $\mathcal{N}(\bm{x};\bm{W},\bm{b})$ could be resolved. That is to say, the model is trained.

Consider the locally supported methods like finite element methods or RKPM, a series of nodes/particles are often defined as the support of approximation. Consider the nodal set $(\bm{x}_1,\bm{x}_2,...,\bm{x}_I,...,\bm{x}_J,...,\bm{x}_{N_\text{P}})$, based on which the coarse-scale approximation function could be recast to the following form
\begin{equation}
\bar{y}^\text{h}(\bm{x}) = \sum_{J = 1}^{N_\text{P}} \Psi_J(\bm{x}, \bm{x}_J) d_J
\end{equation}
where the kernel function $\Psi_J(\bm{x}, \bm{x}_J)$ is compactly defined within the neighborhood domain of $\bm{x}_J$. With proper establishment of kernel functions $\Psi_J(\bm{x}, \bm{x}_J)$, the interpolation is not difficult to obtain as follows
\begin{equation}
y(\bm{x}_I) = y_I = \sum_{J = 1}^{N_\text{P}} \Psi_J(\bm{x}_I, \bm{x}_J) d_J
\end{equation} 
That is to say, the target function is accurately reproduced on the nodal set. If the Kronecker delta property is reached for the kernel function as follows
\begin{equation}
\Psi_J(\bm{x}_I, \bm{x}_J)=
\begin{cases}
1~,~~\bm{x}_I = \bm{x}_J\\
0~,~~\bm{x}_I \neq \bm{x}_J
\end{cases}
\end{equation}
The coefficient will be the value of function at the node, i.e. $d_J = y_J$.

Moreover, we could introduce the additional constraints for the neural network $\mathcal{N}(\bm{x};\bm{W},\bm{b})$ as follows
\begin{equation}\label{eq:constraint_d}
\mathcal{N}(\bm{x}_I;\bm{W},\bm{b}) = 0~,~~\forall \bm{x}_I \in \{\bm{x}_1,\bm{x}_2,...,\bm{x}_I,...,\bm{x}_{N_\text{P}} \}
\end{equation}
Then we will reach a rather good property. That is to say, the fine-scale approximation based on the neural network would not intervene the coarse-scale interpolation, especially for the values of interpolation coefficients $\bm{d} = (d_1, d_2, ....,d_{N_\text{P}})^{\text{T}}$. The NN enhancement full-fill the constrains in Eq.(\ref{eq:constraint_d}) is named as residual free NN, which is the NN version of the residual free bubble (See \cite{HUGHES20054135, Cangiani2007, DOLBOW20083751}). And the residual free property is implemented based on the loss function. In this circumstance, the loss function in Eq. (\ref{eq:multiscale_loss_function}) turns to the residual free version as follows
\begin{equation}
\mathbf{Loss}' =\left\|y(\bm{x}) -\sum_{J = 1}^{N_\text{P}} \Psi_J(\bm{x}) d_J - \mathcal{N}(\bm{x};\bm{W},\bm{b}) \right\|^2_{L^2}
+ \sum_{I = 1}^{N_\text{P}} \beta_\text{d} \left[ \mathcal{N}(\bm{x}_I;\bm{W},\bm{b}) \right]^2
\end{equation}
where $\beta_\text{d}$ is the penalty.

\subsection{Numerical strategy for training}

For the loss function in Eq.(\ref{eq:multiscale_loss_function}), the unknowns  $(\bm{d},\bm{W},\bm{b})$ should be solved by minimizing the loss function. By taking the advantages of the multi-scale feature of the problem, we try to solve different parts of the unknowns in different scales.

By recalling the multi-scale methods in literature, we consider the simple but efficient strategy, e.g. the hierarchical strategy, as the first move. According to the hierarchical strategy, the coupling between scales has been simplified. The procedure is as follows:

\begin{itemize}
	\item Step I: the approximation problem is solved without considering of NN. The loss function is adopted as follows:
	\begin{equation}
		\mathbf{Loss}_{\text{I}} = \|y(\bm{x}) - \bar{y}^\text{h}(\bm{x})\|^2_{L^2}
		=\left\|y(\bm{x}) -\sum_{J = 1}^{N_\text{P}} \Psi_J(\bm{x}) d_J\right\|^2_{L^2}
	\end{equation}
	the values of approximate coefficients $\bm{d} = (d_1, d_2, ....,d_{N_\text{P}})^{\text{T}}$ could be easily solved based on conventional approximation methods.
	
	\item Step II: take the residual 
	\begin{equation}
		\bar{y}^\text{r} = y(\bm{x}) - \bar{y}^\text{h}(\bm{x}) = y(\bm{x}) - \sum_{J = 1}^{N_\text{P}} \Psi_J(\bm{x}) d_J
	\end{equation}
	as the target function. The loss function is adopted as follows:
	\begin{equation}
	\mathbf{Loss}_{\text{II}} = \|\bar{y}^\text{r} - \tilde{y}^\text{h}(\bm{x})\|^2_{L^2}
	=\left\| \bar{y}^\text{r} - \mathcal{N}(\bm{x};\bm{W},\bm{b}) \right\|^2_{L^2}
	\end{equation}
	or the residual free version:
	\begin{equation}\label{eq:res_free_loss_int}
	\begin{split}
	\mathbf{Loss}'_{\text{II}} &= \|\bar{y}^\text{r} - \tilde{y}^\text{h}(\bm{x})\|^2_{L^2} + \sum_{I = 1}^{N_\text{P}} \beta_\text{d} \left[ \mathcal{N}(\bm{x}_I;\bm{W},\bm{b}) \right]^2\\
	& =\left\| \bar{y}^\text{r} - \mathcal{N}(\bm{x};\bm{W},\bm{b}) \right\|^2_{L^2}
	+ \sum_{I = 1}^{N_\text{P}} \beta_\text{d} \left[ \mathcal{N}(\bm{x}_I;\bm{W},\bm{b}) \right]^2
	\end{split}	
	\end{equation}	
	The parameters $(\bm{W},\bm{b})$ within the neural network could be determined based on regular train process.
	
	\item Step III: the hierarchical solution for the multi-scale approximation problem is 
	\begin{equation}
		y(\bm{x}) \approx y^\text{h}(\bm{x}) = \sum_{J = 1}^{N_\text{P}} \Psi_J(\bm{x}) d_J + \mathcal{N}(\bm{x};\bm{W},\bm{b})
	\end{equation}
\end{itemize}

As can be seen, the hierarchical strategy suggests a one-way method without any iterations. Although it is simple enough, the accuracy could fulfill many problems in science and engineering. For the residual free NN enhancement, it is easily to reach the conclusion that no iteration is needed for the approximation problem.

\section{Dirichlet problem in multiscale variational form}
\label{sec:III}

\subsection{Dirichlet problem}

Consider an open bounded domain $\Omega \in \mathbb{R}^d$, for which $d \ge 1$ denotes the dimension of the domain. The boundary of $\Omega$ is denoted by $\Gamma$ or $\partial \Omega$. The boundary value problem is expressed as
\begin{equation}\label{eq:Dirichlet_problem}
\begin{cases}
\mathcal{L} \bm{u} = \bm{f}~~~\text{in}~~~\Omega \\
\bm{u} = \bm{g}~~~\text{on}~~~\Gamma_\text{g}
\end{cases}
\end{equation}
where $\bm{f}$ and $\bm{g}$ are given functions defined in $\Omega$ and $\Gamma_\text{g} \in \Gamma$, respectively. In the present work, we consider the differential operator $\mathcal{L}$ as a second order differential operator.

Consider the trial function space $\mathcal{U}$  and the weighting function space $\mathcal{V}$. The variational form of the BVP (\ref{eq:Dirichlet_problem}) could be expressed as follows
\begin{equation}\label{eq:var_form}
a(\bm{v},\bm{u}) = (\bm{v},\bm{f})
\end{equation}
where the functions $\bm{u} \in \mathcal{U}$ and $\bm{v} \in \mathcal{V}$; $(\cdot,\cdot)$ is the standard inner production in $\mathbb{R}^d$; and the bi-linear form is
\begin{equation}
a(\bm{v},\bm{u}) = (\bm{v},\mathcal{L} \bm{u})
\end{equation}
On the Dirichlet boundary one obtains
\begin{equation}
\begin{cases}
\bm{u} = \bm{g} ~~~\text{on}~~~\Gamma_\text{g}~,~~\forall~~\bm{u} \in \mathcal{U} \\
\bm{v} = \bm{0} ~~~\text{on}~~~\Gamma_\text{g}~,~~\forall~~\bm{v} \in \mathcal{V} \\
\end{cases}
\end{equation}

\subsection{Variational multi-scale form}
According to the celebrated works of Hughes and others \cite{HUGHES19983, GARIKIPATI200039, Cangiani2007, DOLBOW20083751}, one could introduce the decomposition form of solution as follows
\begin{equation}\label{eq:u+u}
\bm{u} = \bar{\bm{u}} + \tilde{\bm{u}}
\end{equation}
where $\bar{\bm{u}}$ and $\tilde{\bm{u}}$ are coarse-scale solution and fine-scale solution, respectively. The weighting function could be also decomposed as follows
\begin{equation}\label{eq:v+v}
\bm{v} = \bar{\bm{v}} + \tilde{\bm{v}}
\end{equation}
The corresponding function spaces are defined as follows
\begin{equation}
\begin{cases}
\mathcal{U} = \bar{\mathcal{U}} \oplus \tilde{\mathcal{U}} \\
\mathcal{V} = \bar{\mathcal{V}} \oplus \tilde{\mathcal{V}}
\end{cases}
\end{equation}
where
\begin{equation}
\begin{cases}
\bar{\bm{u}} \in \bar{\mathcal{U}}~,~~\tilde{\bm{u}} \in \tilde{\mathcal{U}} \\
\bar{\bm{v}} \in \bar{\mathcal{V}}~,~~\tilde{\bm{v}} \in \tilde{\mathcal{V}} 
\end{cases}
\end{equation}
For the Dirichlet boundary we adopt the following forms
\begin{equation}
\begin{cases}
\bar{\bm{u}} = \bm{g}~~~\text{on}~~~\Gamma_\text{g} \\
\tilde{\bm{u}} = \bm{0}~~~\text{on}~~~\Gamma_\text{g} \\
\bar{\bm{v}} = \bm{0}~~~\text{on}~~~\Gamma_\text{g} \\
\tilde{\bm{v}} = \bm{0}~~~\text{on}~~~\Gamma_\text{g}
\end{cases}
\end{equation}

Substitute Eqs. (\ref{eq:u+u}) and (\ref{eq:v+v}) into Eq. (\ref{eq:var_form}), we obain the multiscale variational form as follows
\begin{equation}\label{eq:multi_v_s}
a(\bar{\bm{v}} + \tilde{\bm{v}},\bar{\bm{u}} + \tilde{\bm{u}}) = (\bar{\bm{v}} + \tilde{\bm{v}},\bm{f})
\end{equation}
By using the independence between $\bar{\bm{v}}$ and $\tilde{\bm{v}}$, two problems are obtained
\begin{itemize}
	\item Problem 1: coarse scale problem
	\begin{equation}\label{eq:coarse_scale_prob}
	\begin{cases}
	a(\bar{\bm{v}},\bar{\bm{u}} ) + a(\bar{\bm{v}}, \tilde{\bm{u}}) = (\bar{\bm{v}},\bm{f}) \\
	\bar{\bm{u}} = \bm{g}~,\bar{\bm{v}}=\bm{0}~~\text{on}~~~\Gamma_\text{g}
	\end{cases}	
	\end{equation}
	\item Problem 2: fine scale problem 
	\begin{equation}\label{eq:fine_scale_prob}
	\begin{cases}
	a(\tilde{\bm{v}}, \tilde{\bm{u}}) + a(\tilde{\bm{v}},\bar{\bm{u}}) = ( \tilde{\bm{v}},\bm{f}) \\
	\tilde{\bm{u}} = \bm{0}~,\tilde{\bm{v}}=\bm{0}~~\text{on}~~~\Gamma_\text{g}
	\end{cases}	
	\end{equation}
\end{itemize}
It is observed that the terms $a(\bar{\bm{v}}, \tilde{\bm{u}})$ and $a(\tilde{\bm{v}},\bar{\bm{u}})$ govern the inter-scale coupling. 

\subsection{Solution strategies}

In the conventional variational multi-scale method \cite{HUGHES19983}, the fine-scale problem in Eq.(\ref{eq:fine_scale_prob}) is solved analytically by using the Green function methods in the following form
\begin{equation}\label{eq:Green_Func_fine}
\tilde{\bm{u}} = -\int_{\Omega} g(x,y) (\mathcal{L}\bar{\bm{u}} - \bm{f})(y) \text{d} \Omega
\end{equation}
where $g(x,y)$ is the Green function of the differential operator $\mathcal{L}$. By substutiting the Green function based fine-scale formal solution in Eq.(\ref{eq:Green_Func_fine}) to the coarse-scale problem in Eq.(\ref{eq:coarse_scale_prob}), the whole problem could be recast to
\begin{equation}\label{eq:VMM_solution}
\begin{cases}
	a(\bar{\bm{v}},\bar{\bm{u}} ) - a(\bar{\bm{v}}, \int_{\Omega} g(x,y) (\mathcal{L}\bar{\bm{u}} - \bm{f})(y) \text{d} \Omega) = (\bar{\bm{v}},\bm{f}) \\
	\bar{\bm{u}} = \bm{g}~,\bar{\bm{v}}=\bm{0}~~\text{on}~~~\Gamma_\text{g}
\end{cases}	
\end{equation}
As we can see, Eq.(\ref{eq:VMM_solution}) is an exact equation for the coarse scale problem with the effects from the fine scale in the non-local form.

In the present work, we shall propose an alternate strategy for variational multi-scale problem. For the first step, we are trying to take the advantage of the conventional numerical methods, such as finite difference methods, finite element methods, meshfree methods and so on. By applying proper conventional methods, we solve the corase scale problem in Eq. (\ref{eq:coarse_scale_prob}) by neglecting the fine scale terms as follows
\begin{equation}\label{eq:coarse_scale_prob_wo_fine}
\begin{cases}
a(\bar{\bm{v}},\bar{\bm{u}} ) = (\bar{\bm{v}},\bm{f}) \\
\bar{\bm{u}} = \bm{g}~,\bar{\bm{v}}=\bm{0}~~\text{on}~~~\Gamma_\text{g}
\end{cases}	
\end{equation}
Actually, it is the conventional numerical solution for the original Dirichlet problem in Eq.(\ref{eq:Dirichlet_problem}).

For the second step, the fine scale problem in Eq.(\ref{eq:fine_scale_prob}) could be solved with the coarse scale term $\bar{\bm{u}}$ given numerically. Here a kind of down-scaling technique may be applied. In the present work, we consider the neural network based solution techniques for the numerical solution in fine-scale. 

Finally, the multi-scale solution could be easily obtained in the form of $\bm{u} = \bar{\bm{u}} + \tilde{\bm{u}}$ (also shown in Eq. (\ref{eq:u+u})). On the other side, the governing PDE of the fine scale problem in Eq.(\ref{eq:fine_scale_prob}) could be recast to
\begin{equation}
a(\tilde{\bm{v}}, \tilde{\bm{u}}) + a(\tilde{\bm{v}},\bar{\bm{u}}) 
= a(\tilde{\bm{v}}, \bar{\bm{u}} + \tilde{\bm{u}}) 
= a(\tilde{\bm{v}}, \bm{u}) 
= ( \tilde{\bm{v}},\bm{f})
\end{equation}
Due to the arbitrariness of $\tilde{\bm{v}}$, the multi-stale solution satisfies the Galerkin form of the original PDE of Dirichlet problem in Eq.(\ref{eq:Dirichlet_problem}). And the essential boundary conditions
\begin{equation}
\bar{\bm{u}} + \tilde{\bm{u}} = \bm{u} = \bm{g}~~~\text{on}~~~\Gamma_\text{g}
\end{equation}
also fulfills. The proposed strategy is decoupled and rather easily to be implemented. In addition, the total solution is of good accuracy in the sense of fine-scale Galerkin form.

\section{Neural network enhanced methods}
\label{sec:IV}

\subsection{NN based enhancement}

The coarse-scale functions are also approximated based on conventional methods in the following form
\begin{equation}
\begin{cases}
\bar{\bm{u}}(\bm{x}) \approx \bar{\bm{u}}^\text{h}(\bm{x}) = \sum_{J = 1}^{N_\text{P}} \Psi_J(\bm{x},\bm{x}_J) \bm{u}_J \\
\bar{\bm{v}}(\bm{x}) \approx \bar{\bm{v}}^\text{h}(\bm{x}) = \sum_{J = 1}^{N_\text{P}} \Psi_J(\bm{x},\bm{x}_J) \bm{v}_J
\end{cases}
\end{equation}
where $\bm{u}_J$ and $\bm{v}_J$ are interpolation coefficients at $\bm{x}_J$. The derivatives with respect to $\bm{x}$ are in the form
\begin{equation}
\begin{cases}
\frac{\partial}{\partial \bm{x}} \bar{\bm{u}}(\bm{x}) \approx \frac{\partial}{\partial \bm{x}} \bar{\bm{u}}^\text{h}(\bm{x}) = \sum_{J = 1}^{N_\text{P}} \frac{\partial}{\partial \bm{x}} \Psi_J(\bm{x},\bm{x}_J) \bm{u}_J \\
\frac{\partial}{\partial \bm{x}} \bar{\bm{v}}(\bm{x}) \approx \frac{\partial}{\partial \bm{x}} \bar{\bm{v}}^\text{h}(\bm{x}) = \sum_{J = 1}^{N_\text{P}} \frac{\partial}{\partial \bm{x}} \Psi_J(\bm{x},\bm{x}_J) \bm{v}_J
\end{cases}
\end{equation}

The fine-scale functions are approximated based on neural networks. Thus we have
\begin{equation}\label{eq:NN_form_fine}
\begin{cases}
\tilde{\bm{u}}(\bm{x}) \approx \tilde{\bm{u}}^\text{h}(\bm{x}) = \mathcal{N}(\bm{x};\bm{W},\bm{b}) \\
\tilde{\bm{v}}(\bm{x}) \approx \tilde{\bm{v}}^\text{h}(\bm{x}) = \mathcal{N}(\bm{x};\bm{W},\bm{b})
\end{cases}
\end{equation}
The derivatives  with respect to $\bm{x}$ are in the form
\begin{equation}
\begin{cases}
\frac{\partial}{\partial \bm{x}} \tilde{\bm{u}}(\bm{x}) \approx \frac{\partial}{\partial \bm{x}} \tilde{\bm{u}}^\text{h}(\bm{x}) = \frac{\partial}{\partial \bm{x}} \mathcal{N}(\bm{x};\bm{W},\bm{b}) \\
\frac{\partial}{\partial \bm{x}} \tilde{\bm{v}}(\bm{x}) \approx \frac{\partial}{\partial \bm{x}} \tilde{\bm{v}}^\text{h}(\bm{x}) = \frac{\partial}{\partial \bm{x}} \mathcal{N}(\bm{x};\bm{W},\bm{b})
\end{cases}
\end{equation}
The constitution of the neural networks $\mathcal{N}(\bm{x};\bm{W},\bm{b})$ could be found in Eqs. (\ref{eq:NN_module}) and (\ref{eq:NN_layer}). Moreover, the automatic differentiation\cite{Automatic_differentiation_2015} could be easily performed for $\mathcal{N}(\bm{x};\bm{W},\bm{b})$ based on the package like TensorFlow\cite{TensorFlow2021} or PyTorch\cite{PyTorch2017}. In the present work, we follow the standard Galerkin method so that the same functional space is used for the solution function space and the test function space.

\subsection{Loss functions in variational forms}

Consider the variational form in Eq.(\ref{eq:fine_scale_prob}). If the different operator corresponding to the bilinear form is self-adjoint, the bilinear form is symmetric and the energy form could be established. We have
\begin{equation}
\Pi (\tilde{\bm{u}}) = \frac{1}{2} \| \tilde{\bm{u}} \|^2_E + a(\bar{\bm{u}},\tilde{\bm{u}}) - ( \tilde{\bm{u}},\bm{f})
\end{equation}
where the energy norm is defined as
\begin{equation}
\|\bm{u}\|^2_E = a(\bm{u}, \bm{u})
\end{equation}
The weak form of the problem recast to
\begin{equation}
\begin{split}
\text{find}~~&\min~\Pi (\tilde{\bm{u}}) \\
\text{s.t.}~~&~\tilde{\bm{u}} = \bm{0}~~\text{on}~~~\Gamma_\text{g}
\end{split}
\end{equation}
The loss function could be given in the following form
\begin{equation}\label{eq:energy_loss}
\mathbf{Loss}_1 = \Pi (\tilde{\bm{u}}) + \frac{\alpha_\text{p}}{2} 
\int_{\Gamma_\text{g}} \tilde{\bm{u}} \cdot \tilde{\bm{u}} \text{d} \Gamma
\end{equation}
where $\alpha_\text{p}$ is a penalty to impose the essential boundary condition. The energy based loss function in Eq.(\ref{eq:energy_loss}) is a good candidate to train the neural networks $\mathcal{N}(\bm{x};\bm{W},\bm{b})$ for approximating the field variable $\tilde{\bm{u}}$, as shown in Eq.(\ref{eq:NN_form_fine}). 

In certain circumstances, external constraints are introduced to the problem. Consider the constraints in the form of Eq. (\ref{eq:constraint_d}), as an example. Then the loss function would be recast to the following form:
\begin{equation}\label{eq:energy_loss_residual_free}
\mathbf{Loss}'_1 = \Pi (\tilde{\bm{u}}) + \frac{\alpha_\text{p}}{2} 
\int_{\Gamma_\text{g}} \tilde{\bm{u}} \cdot \tilde{\bm{u}} \text{d} \Gamma 
+ \sum_{I = 1}^{N_\text{P}} \beta_\text{d} \left[ \mathcal{N}(\bm{x}_I;\bm{W},\bm{b}) \right]^2
\end{equation}
The constraints in other forms could be implemented in the similar way.

If the differential operator associate with the bilinear form is not self-adjoint, the bilinear form is not symmetric. The energy form might not exist for this case. We reform Eq.(\ref{eq:fine_scale_prob}) into its strong form as follows
\begin{equation}\label{eq:strong_form}
\begin{cases}
\mathcal{F}(\tilde{\bm{u}}) =
\mathcal{L} \tilde{\bm{u}} - \mathcal{L} \bar{\bm{u}} -\bm{f} = \bm{0} \\
\tilde{\bm{u}} = \bm{0}~~\text{on}~~~\Gamma_\text{g}
\end{cases}	
\end{equation}
By referring to PINN \cite{RAISSI2019686}, the loss function could be developed based on collocation methods. We have
\begin{equation}
\mathbf{Loss}_2 = \frac{1}{N_\text{f}} \sum_{i = 1}^{N_\text{f}} \left[ \mathcal{F}(\tilde{\bm{u}}_i) \right]^2 +
\frac{\beta_\text{c}}{N_\text{u}}\sum_{j = 1}^{N_\text{u}} \left[ \tilde{\bm{u}}_j \right]^2
\end{equation}
Here we introduce a coefficient $\beta_\text{c}$ to balance the errors between inside the domain and on the boundary. Similarly, the version of loss function with external constraints could be expressed as follows
\begin{equation}
\mathbf{Loss}'_2 = \frac{1}{N_\text{f}} \sum_{i = 1}^{N_\text{f}} \left[ \mathcal{F}(\tilde{\bm{u}}_i) \right]^2 +
\frac{\beta_\text{c}}{N_\text{u}}\sum_{j = 1}^{N_\text{u}} \left[ \tilde{\bm{u}}_j \right]^2
+ \frac{\beta_\text{d}}{N_\text{P}}\sum_{I = 1}^{N_\text{P}} \left[ \mathcal{N}(\bm{x}_I;\bm{W},\bm{b}) \right]^2
\end{equation}

The neural networks and the loss functions could be implemented by popular machine learning packages, e.g. TensorFlow\cite{TensorFlow2021} or PyTorch\cite{PyTorch2017}, with automatic differentiation\cite{Automatic_differentiation_2015}.

\subsection{Implementation and algorithm set-up}

We still consider the hierarchical scheme shown in the following items:
\begin{itemize}
	\item Step I: to solve Problem 1 subjected to $\tilde{\bm{u}} \leftarrow \tilde{\bm{u}}^\text{h}_0 = \bm{0}$ with conventional numerical methods, e.g. finite element method, meshfree methods, and so on.
	\begin{equation}\label{eq:coarse_alg}
	\left.
	\begin{split}
	&a(\bar{\bm{v}},\bar{\bm{u}} ) = (\bar{\bm{v}},\bm{f})~~\text{in}~~~\Omega \\
	&\bar{\bm{u}} = \bm{g}~,\bar{\bm{v}}=\bm{0}~~\text{on}~~~\Gamma_\text{g}
	\end{split}
	\right\rbrace
	\Longrightarrow
	\bar{\bm{u}}^\text{h}_1 \leftarrow \bar{\bm{u}} = \sum_{J = 1}^{N_\text{P}} \Psi_J(\bm{x},\bm{x}_J) \bm{u}_J
	\end{equation}
	
	\item Step II: to solve Problem 2 subjected to $\bar{\bm{u}} = \bar{\bm{u}}^\text{h}_1$ with NN based approximations.
	\begin{equation}\label{eq:fine_alg}
	\left.
	\begin{split}
	a(\tilde{\bm{v}}, \tilde{\bm{u}}) + a(\tilde{\bm{v}},\bar{\bm{u}}^\text{h}_1) = ( \tilde{\bm{v}},\bm{f})~~\text{in}~~\Omega_{\text{p}} \\
	\tilde{\bm{u}} = \bm{0}~,\tilde{\bm{v}}=\bm{0}~~\text{on}~~~\Gamma_\text{g}
	\end{split}
	\right\rbrace
	\Longrightarrow
	\tilde{\bm{u}}^\text{h}_1 \leftarrow \tilde{\bm{u}} = \mathcal{N}(\bm{x}_I;\bm{W},\bm{b})
	\end{equation}
	
	\item Step III: the hierarchical solution for the multi-scale approximation problem is 
	\begin{equation}\label{eq:total_solution_alg}
	\bm{u} = \bar{\bm{u}} + \tilde{\bm{u}}
	\approx \bar{\bm{u}}^\text{h}_1 + \tilde{\bm{u}}^\text{h}_1  = \sum_{J = 1}^{N_\text{P}} \Psi_J(\bm{x},\bm{x}_J) \bm{u}_J + \mathcal{N}(\bm{x}_I;\bm{W},\bm{b})
	\end{equation}

\end{itemize}

Although the hierarchical scheme is rather simple, it is suitable for many problems and of good accuracy. In the present work, we just focus on the hierarchical scheme. Also, the staggered methods with iterations between coarse scale problem and fine scale problem may be a step-forward and will be investigated in the future work. In this circumstance, attentions should be paid to the term of $a(\bar{\bm{v}}, \tilde{\bm{u}})$, which carries the influence from the fine scale to the coarse scale.

\section{Case studies}
\label{sec:V}

\subsection{Approximation problems}

\subsubsection{Function approximation in 1D}

Firstly, consider the continuous target function in the following form:
\begin{equation}
f(x) = 0.5 x + x^3 + \tanh (10x)~,~~x \in [-1,1]
\end{equation}

\begin{figure}[htbp]
	\centering
	\includegraphics[width=0.7\textwidth]{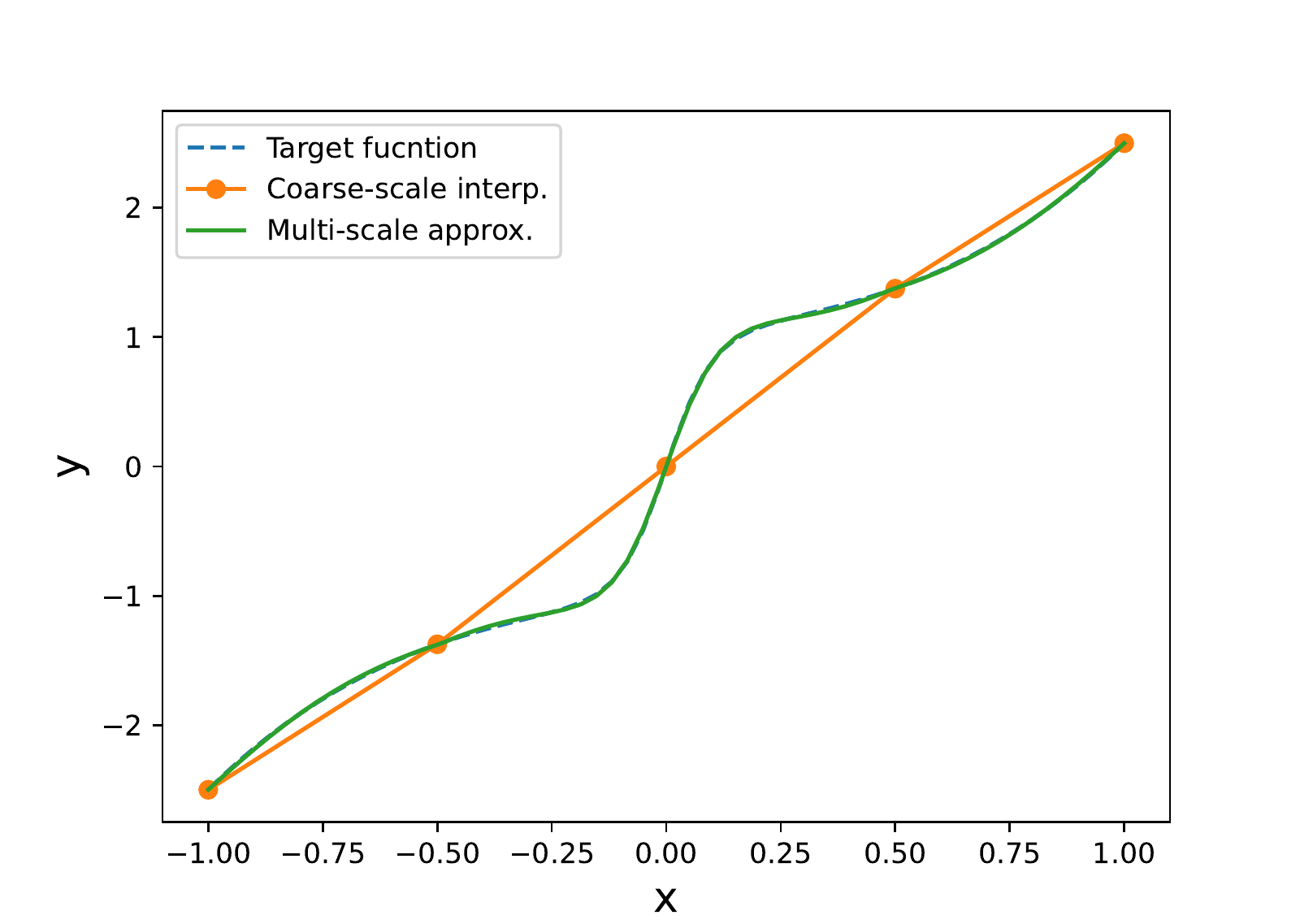}
	\caption{Continuous function approximation based on multiscale-NN framework}
	\label{fig:1}
\end{figure}

\begin{figure}[htbp]
	\centering
	\includegraphics[width=0.7\textwidth]{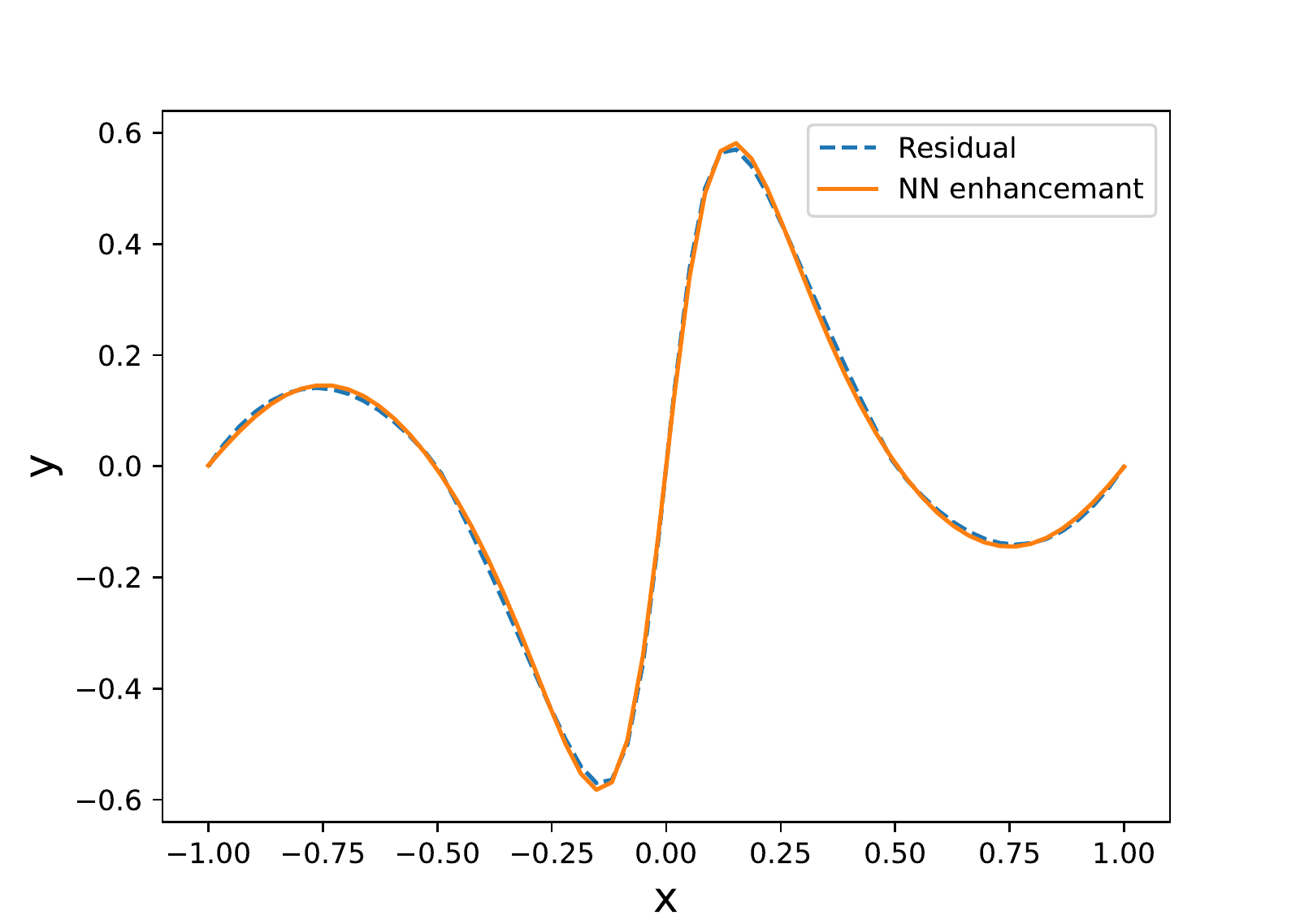}
	\caption{Residual after polynomial approximation}
	\label{fig:2}
\end{figure}

We approximate the target function based on the multiscale-NN framework. In the coarse scale, the domain $[-1,1]$ is evenly divided into four elements. In each element, the linear shape function is applied. The piece-wise linear approximation of the target function could be seen in Fig.\ref{fig:1}. Because we only use four elements for the coarse scale approximation, the residual error is conspicuous. In the fine scale, the residual error from the coarse scale is approximated based on neural networks. A neural network with 2 hidden dense layers is developed for the fine-scale approximation. The numbers of neurons in the hidden layers are $(4,7)$ and the total number of parameters is 51. The activation functions are adopted as Sigmoid for the hidden layers and linear for the output layer. The residual free form of Loss function expressed by Eq.(\ref{eq:res_free_loss_int}) is adopted. And 60 integration points are equally distributed within the domain $[-1,1]$ for the calculation of the loss function based on nodal integration (see 
Appendix). Trained by the Adam trainer for 18000 epochs, we obtain the $L_2$ error in the level of $\text{O}(10^{-5})$. As can be seen in Fig.\ref{fig:2}, the residual between the target function and the coarse scale approximation has been well represented. The superposition of the coarse scale result and the fine scales result yields the multi-scale results, as shown in Fig.\ref{fig:1}. It has no doubt that the multi-scale result is of good accuracy.

Second, consider discontinuous target function expressed as follows:
\begin{equation}\label{eq:dis-cont_func}
f(x) = \frac{x - 1 }{2} + \text{H}(x) ~,~~x \in [-1,1]
\end{equation}
where $\text{H}(\cdot)$ denotes the Heaviside function, which is a jump of magnitude 1 at $x=0$.

\begin{figure}[htbp]
	\centering
	\includegraphics[width=0.7\textwidth]{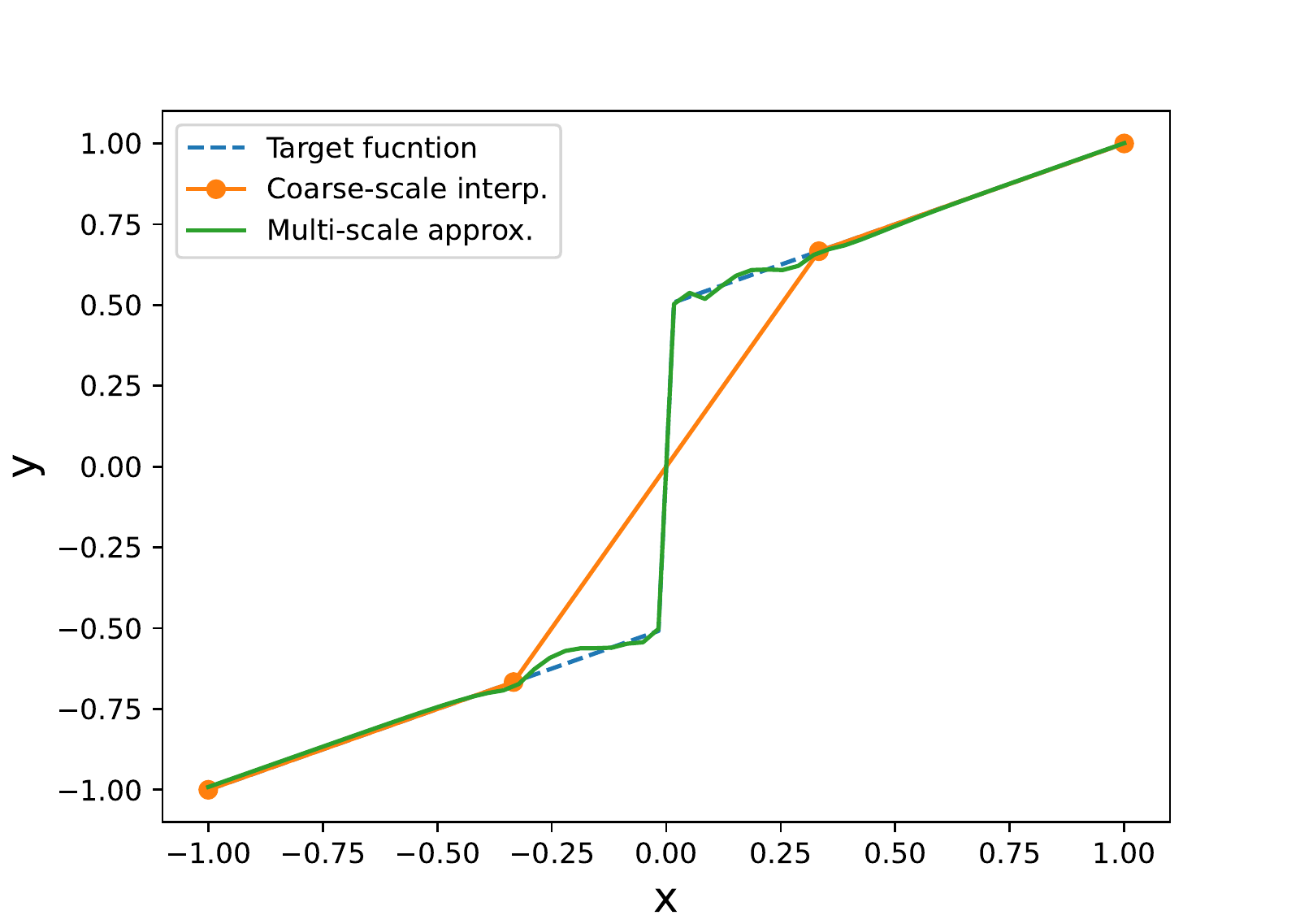}
	\caption{Dis-continuous function approximation based on multiscale-NN framework}
	\label{fig:3}
\end{figure}

\begin{figure}[htbp]
	\centering
	\includegraphics[width=0.7\textwidth]{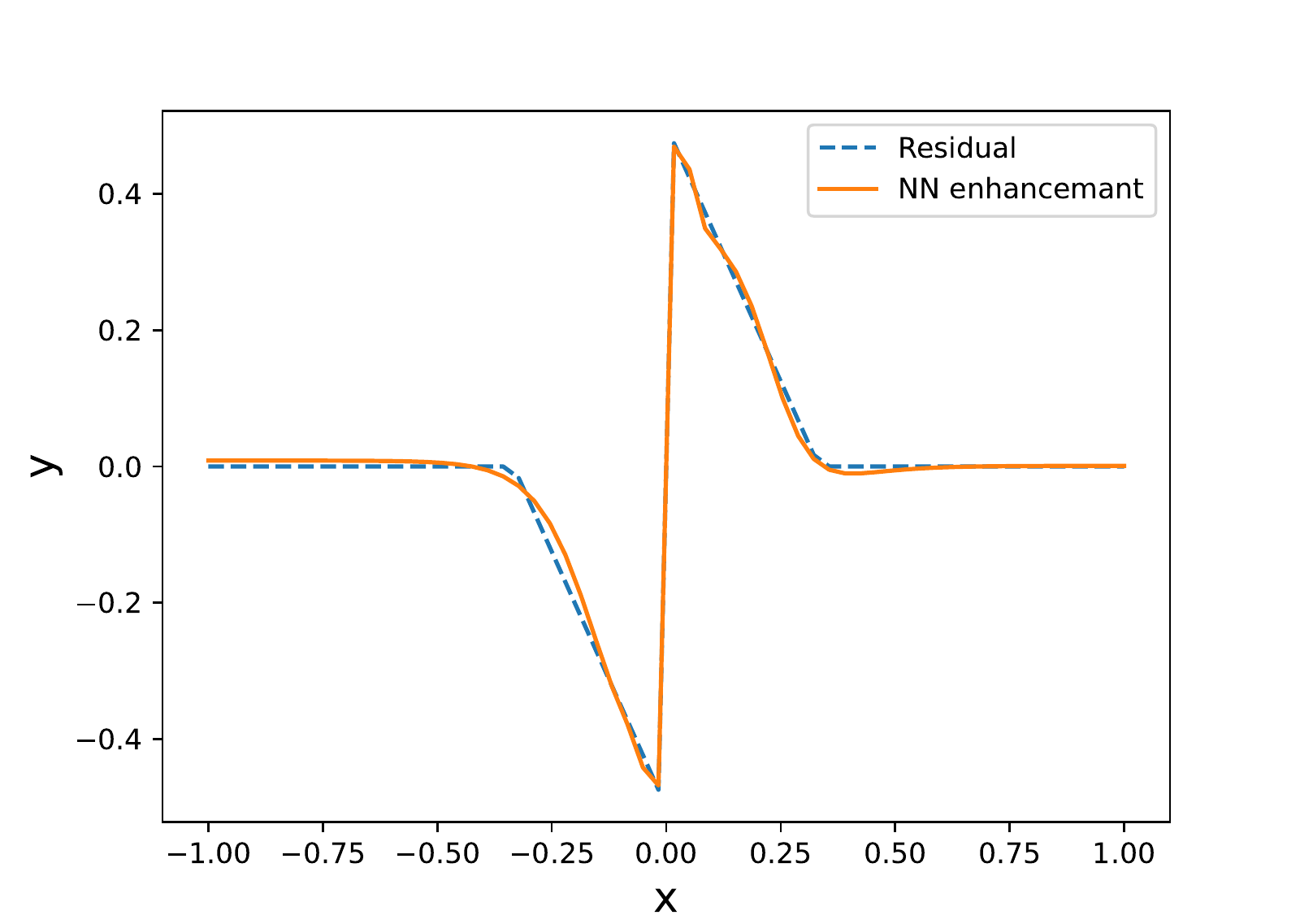}
	\caption{Residual after polynomial approximation}
	\label{fig:4}
\end{figure}

Apply the multiscale framework for the approximation of function in Eq.(\ref{eq:dis-cont_func}). We use three linear elements through $[-1,1]$ as coarse scale approximation, which is very coarse. As can be seen in Fig.\ref{fig:3}, the elements on the left and right give accurate approximation because the target functions are linear within the domain of the corresponding elements. In the middle element, the jump cannot be well captured by the linear shape function. Prominent residual errors between the target function and the linear approximation are clearly observed. Based on the multiscale framework proposed in the present work, the residual between the target function and the coarse scale results is approximated by the NN enhancement in fine scale. A neural network with two hidden layers is developed with the numbers of neurons $(4,8)$. The total number of parameters is 57. The activation functions are adopted as Sigmoid for the first and second hidden layers. The residual free form of Loss function expressed by Eq.(\ref{eq:res_free_loss_int}) is also adopted. And totally 60 integration points are equally distributed within the domain $[-1,1]$. The Adam trainer is used and the total number of epochs is set to be 18000. As can be seen in Fig.\ref{fig:4}, the residual could be well captured by the NN enhancement with the $L_2$ error in the level of $\text{O}(10^{-3})$. As shown in Fig. \ref{fig:3}, the target function is well approximated by the multiscale superposition of coarse scale results and NN enhancement.

\subsubsection{Function approximation in 2D}

\begin{figure}[htbp]
	\centering
	\subfigure[Target function]{
		\label{fig:2D_cont_target_surf} 
		\includegraphics[width=0.48\textwidth]{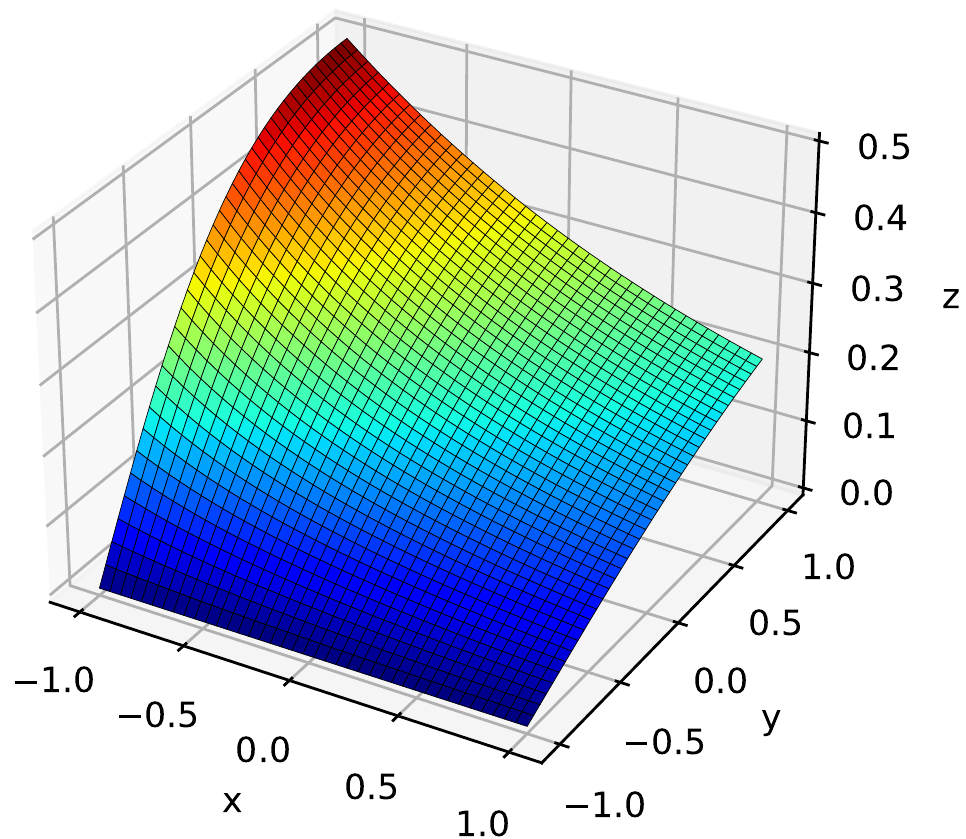}}
	\subfigure[Coarse-scale interpolation]{
		\label{fig:2D_cont_coarsescale} 
		\includegraphics[width=0.48\textwidth]{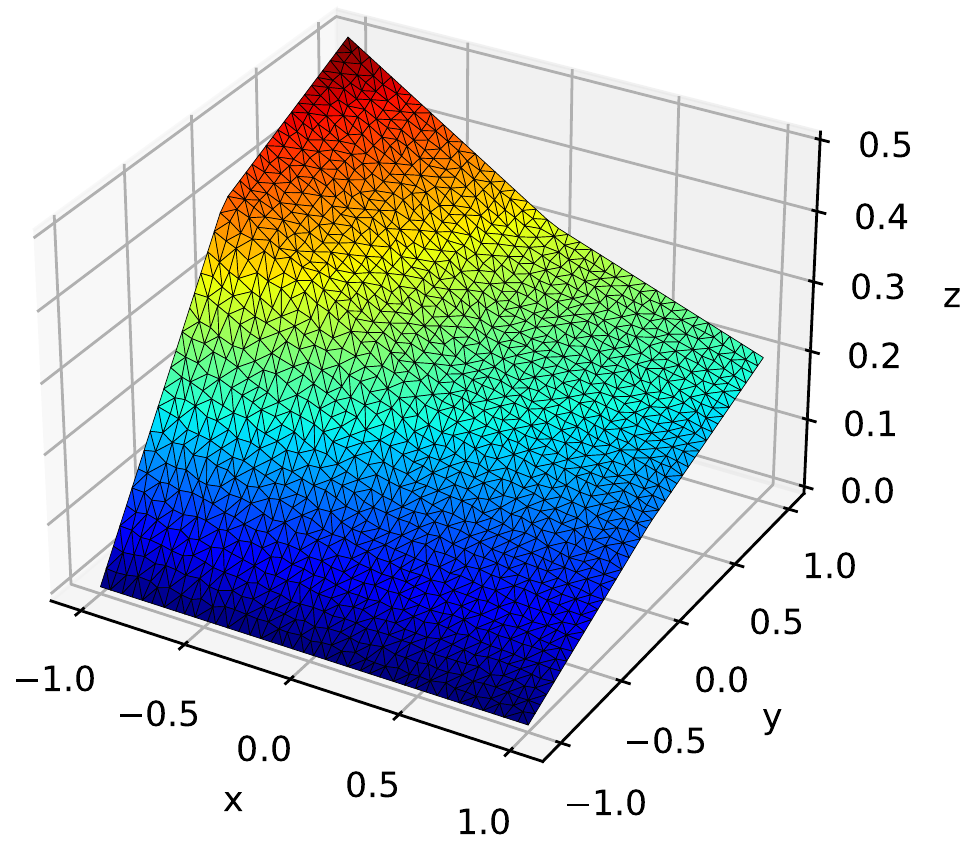}}
	\subfigure[Multi-scale approximation]{
		\label{fig:2D_cont_multiscale} 
		\includegraphics[width=0.48\textwidth]{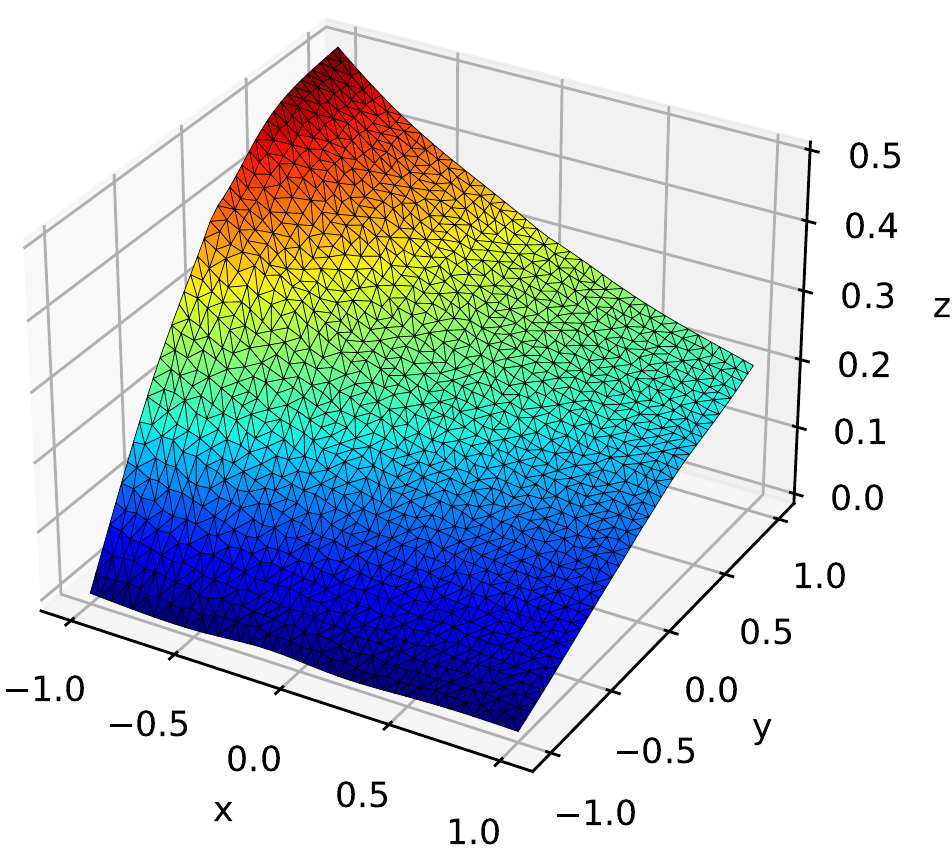}}
	\caption{\label{fig:2D_cont_target}2D continuous function}
\end{figure}

Consider the 2D continuous target function as follows
\begin{equation}
f(x,y) = \frac{2(1+y)}{(3+x)^2 + (1+y)^2}
~,~~ (x,y) \in [-1,1] \otimes [-1,1]
\end{equation}
Its 3D plot is shown in Fig. \ref{fig:2D_cont_target_surf}.

\begin{figure}[htbp]
	\centering
	\subfigure[Residual]{
		\label{fig:2D_cont_rsd} 
		\includegraphics[width=0.48\textwidth]{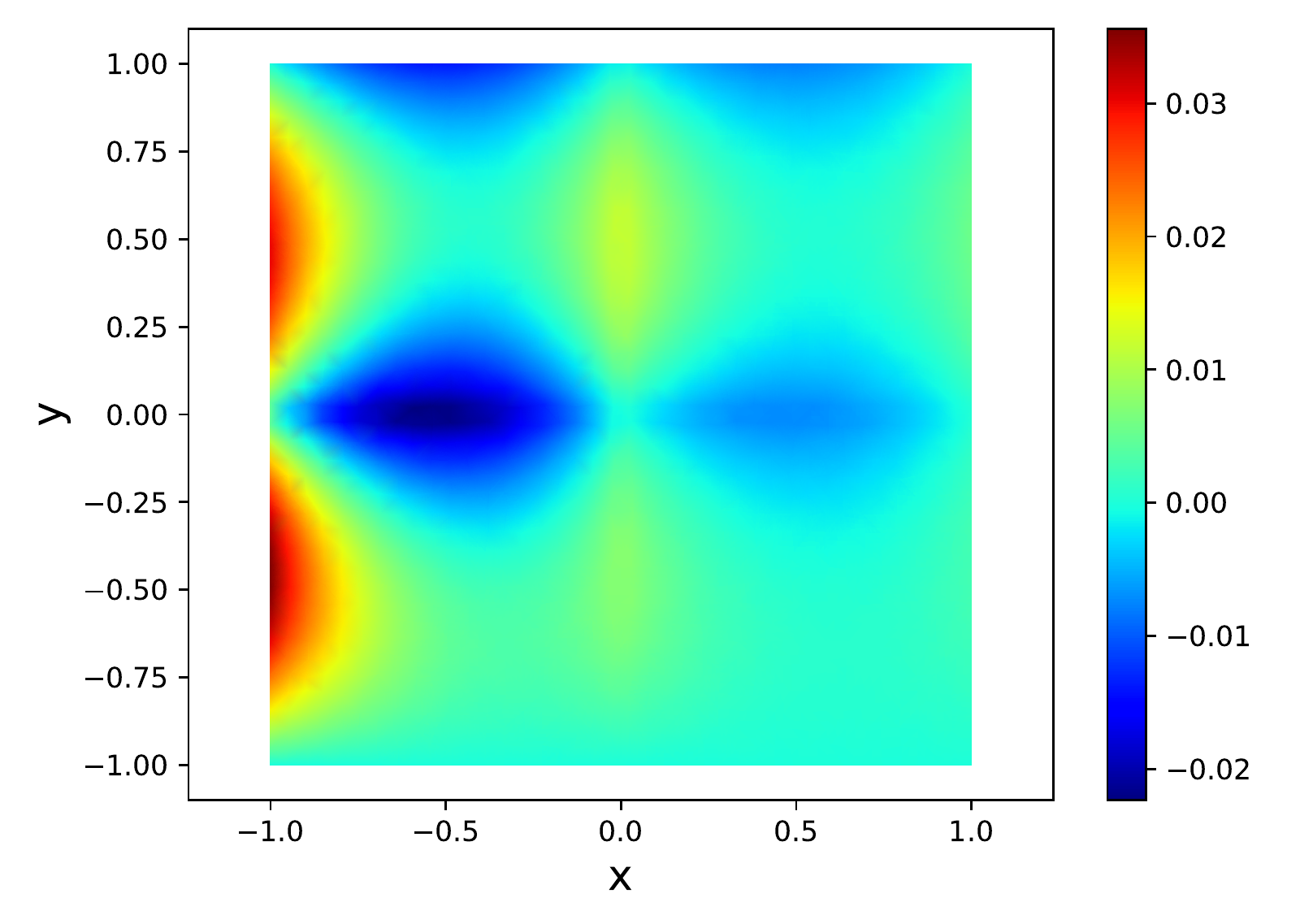}}
	\subfigure[NN enhancement]{
		\label{fig:2D_cont_NN} 
		\includegraphics[width=0.48\textwidth]{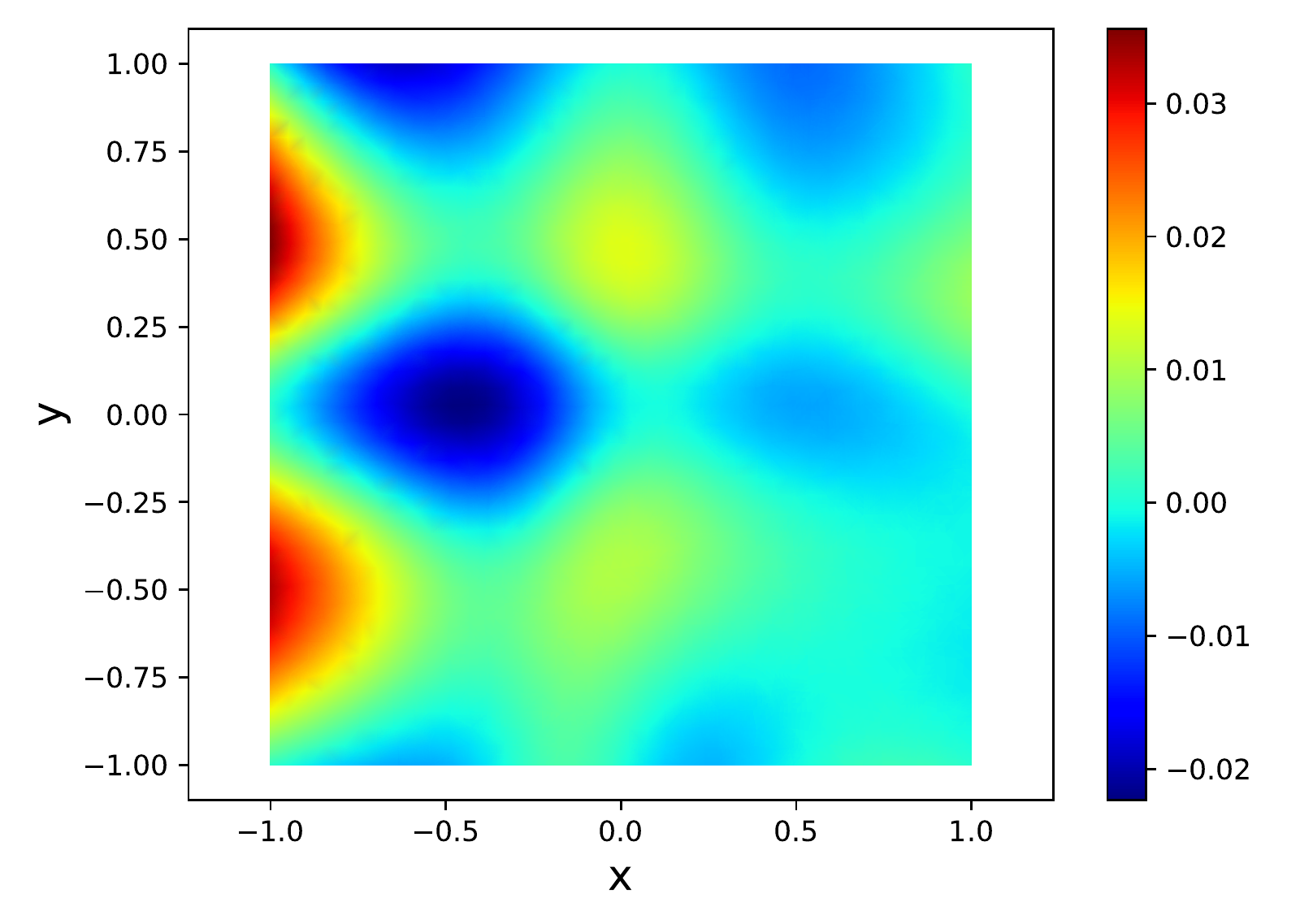}}
	\caption{\label{fig:2D_cont_residual} Residual between target function and coarse scale results}
\end{figure}

Totally 2$\times$2 elements are generated for the coarse scale approximation. For each element, the shape functions are established based on the four corner nodes. It is the standard four-node Lagrange element for the analysis in 2D. The result of coarse scale approximation based on 2$\times$2 four-node-elements could be found in  Fig.\ref{fig:2D_cont_coarsescale}. Comparing Fig.\ref{fig:2D_cont_target_surf} and Fig.\ref{fig:2D_cont_coarsescale}, the difference between them could be clearly observed. The residual between the target function and the coarse-scale approximation is shown in Fig.\ref{fig:2D_cont_rsd}. According to the multiscale framework proposed in the present paper, the residual is approximated by a neural network as the enhancement. Three dense hidden layers are developed with $(8,18,5)$ neurons, respectively. The simple Sigmoid activation function is adopted for each neuron in the hidden layers. The total number of parameters is 205. Still, the residual free loss function in the form of Eq.(\ref{eq:res_free_loss_int}) is used. We employ 40$\times$40 integration points within the domain $[-1,1] \otimes [-1,1]$ for the calculation of loss function by using nodal integration (see 
Appendix). It has no doubt that the residual could be well captured by the NN enhancement, as shown in Fig.\ref{fig:2D_cont_residual}. After being trained by the Adam trainer for 50000 epochs, the $L_2$ error of the NN is reduced to the level of $\text{O}(10^{-6})$. Figure \ref{fig:2D_cont_multiscale} depicts the results of multi-scale approximation, e.g. the superposition of coarse scale results and NN enhancement. It is of reasonable accuracy and agrees well with the target function in Fig.\ref{fig:2D_cont_target_surf}.

\begin{figure}[htbp]
	\centering
	\subfigure[Target function]{
		\label{fig:2D_discont_target_surf} 
		\includegraphics[width=0.48\textwidth]{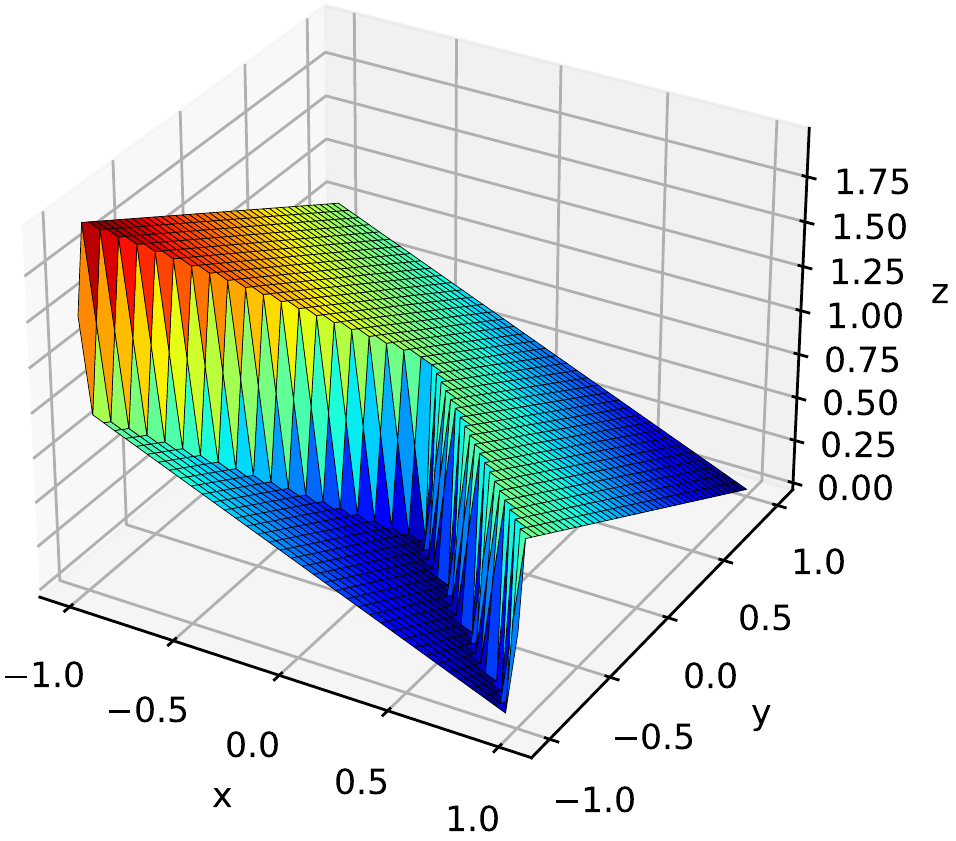}}
	\subfigure[Coarse-scale interpolation]{
		\label{fig:2D_discont_coarsescale} 
		\includegraphics[width=0.48\textwidth]{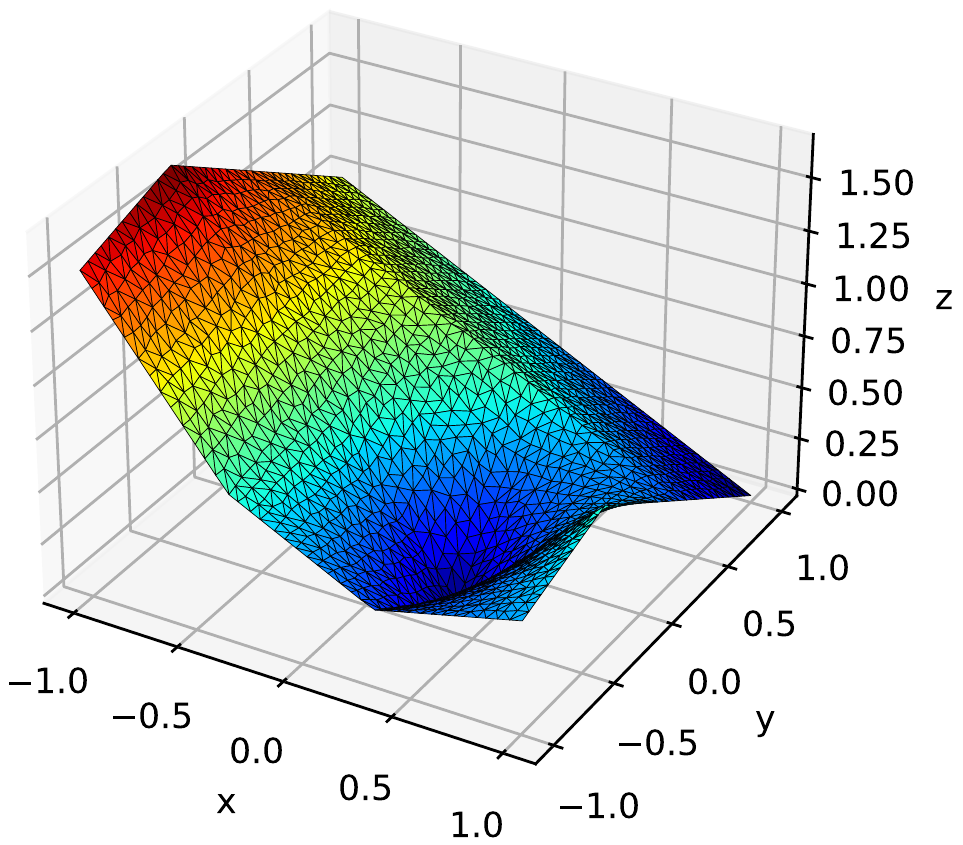}}
	\subfigure[Multi-scale approximation]{
		\label{fig:2D_discont_multiscale} 
		\includegraphics[width=0.48\textwidth]{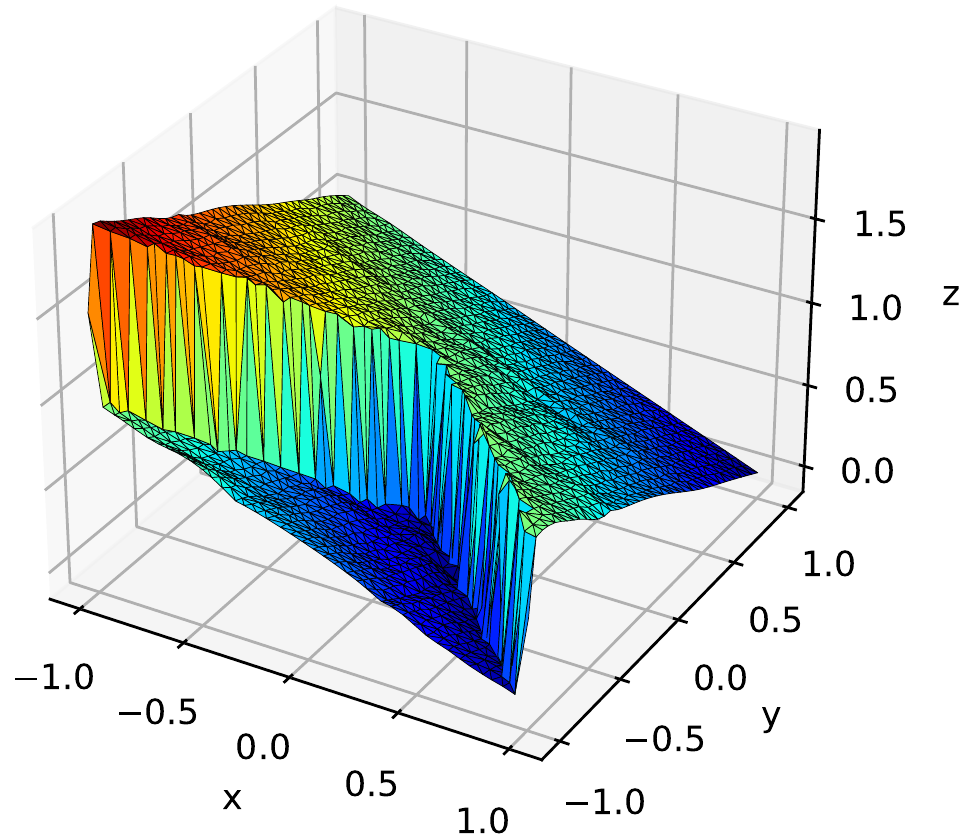}}
	\caption{\label{fig:2D_discont_target}2D discontinuous function}
\end{figure}

Moreover, we consider the 2D discontinuous target function as follows
\begin{equation}
f(x,y) = \text{H}(|x|+y) - \frac{x+y}{2}
~,~~ (x,y) \in [-1,1] \otimes [-1,1]
\end{equation}
As can be seen in Fig. \ref{fig:2D_discont_target_surf}, this target function is strongly discontinuous. 

\begin{figure}[htbp]
	\centering
	\subfigure[Residual]{
		\label{fig:2D_discont_rsd} 
		\includegraphics[width=0.48\textwidth]{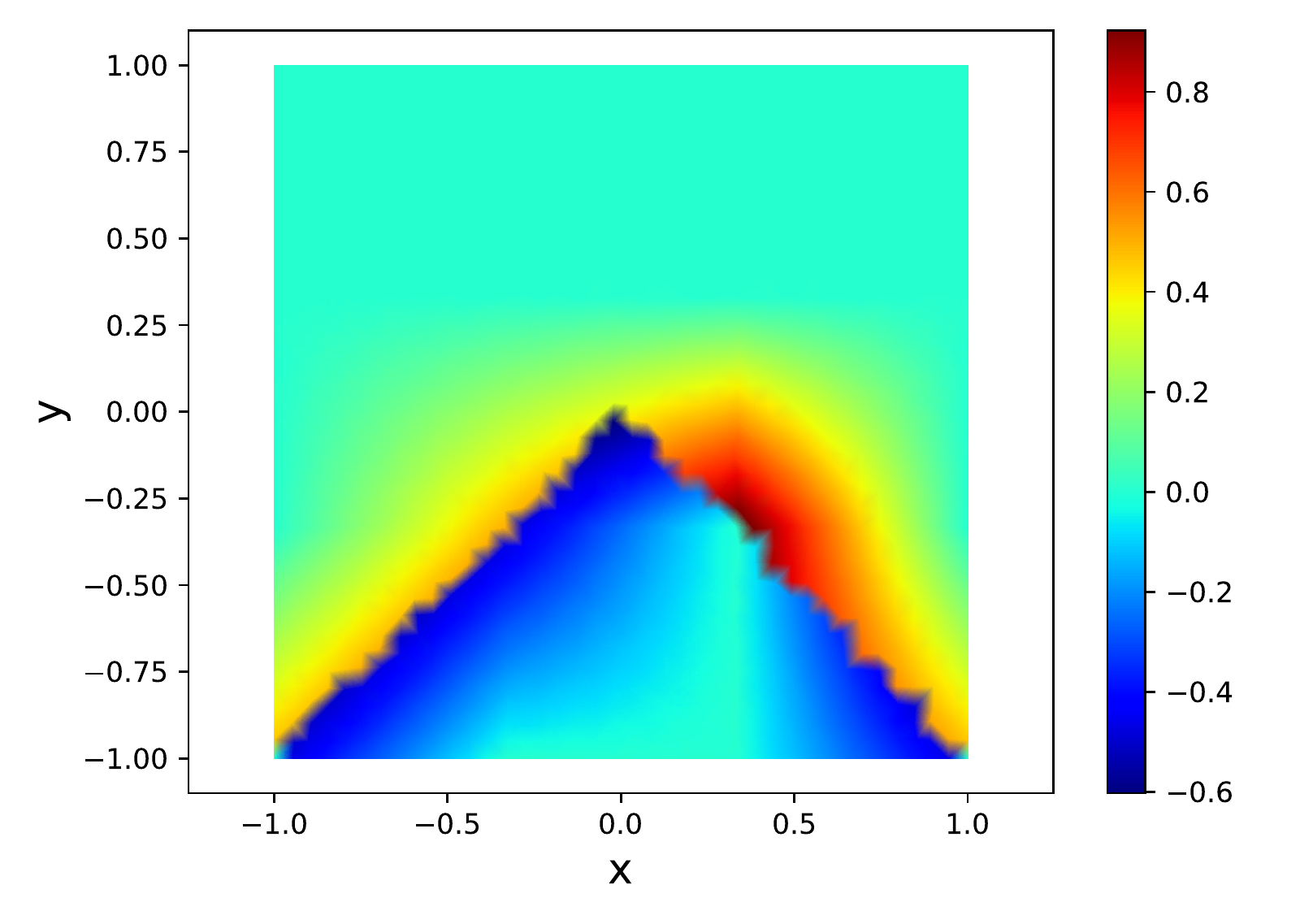}}
	\subfigure[NN enhancement]{
		\label{fig:2D_discont_NN} 
		\includegraphics[width=0.48\textwidth]{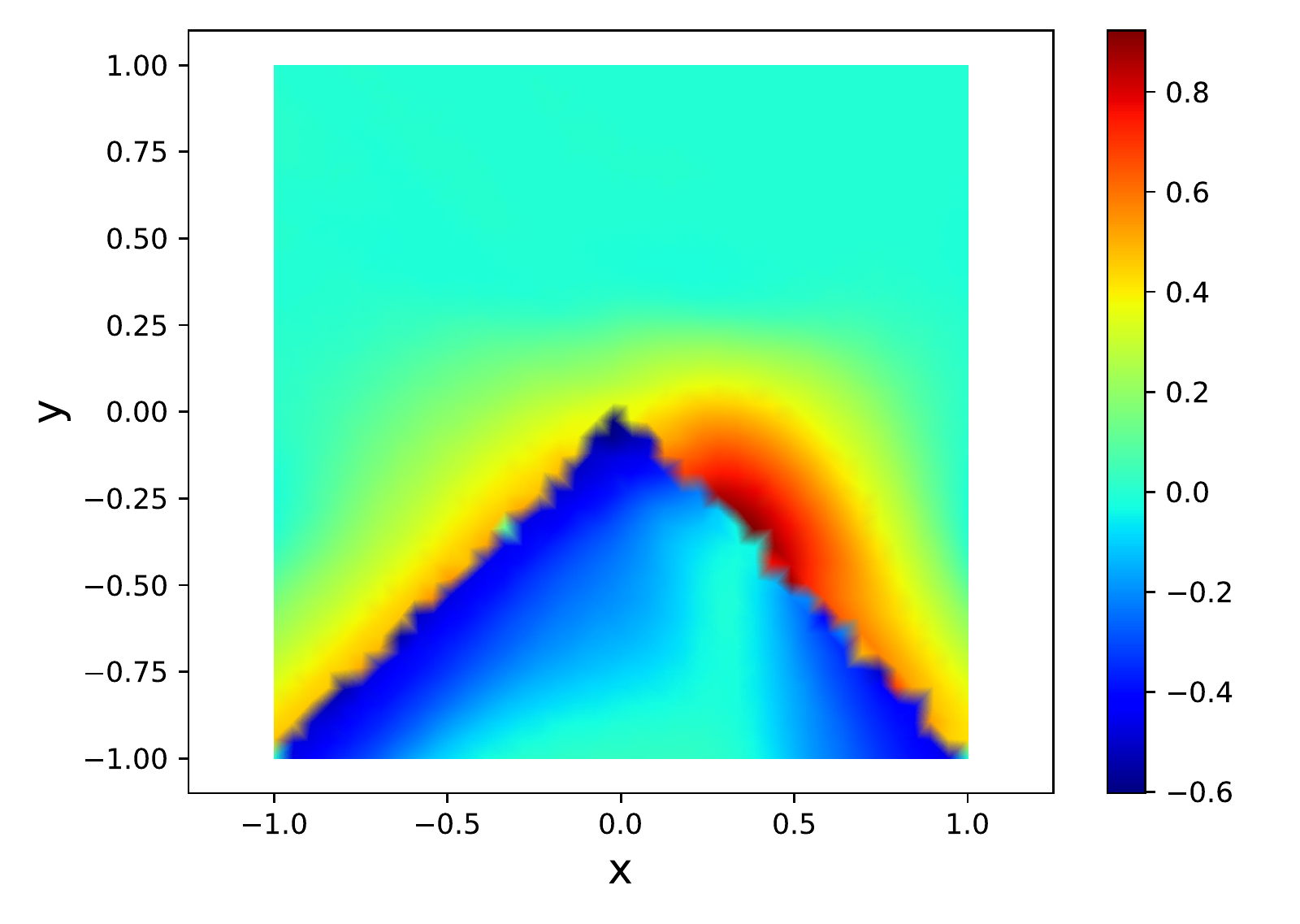}}
	\caption{\label{fig:2D_discont_residual} Residual between target function and coarse scale results}
\end{figure}

For the first step, we develop 3$\times$3 four-node elements for the coarse scale approximation. As shown in Fig. \ref{fig:2D_discont_coarsescale}, the result of the coarse scale approximation is far from the target function in Fig. \ref{fig:2D_discont_target_surf}. The discontinuity can hardly be captured by the regular finite element with continuous shape functions. The residual between the exact target function and the coarse scale approximation is shown in Fig. \ref{fig:2D_discont_rsd}. According to the proposed multiscale framework, the residual could be reproduced by NN enhancement. A neural network with four dense hidden layers is developed. The neurons in the corresponding hidden layers are $(8,25,10,5)$ and the simple Sigmoid activation function is setup for each neuron in the hidden layers. The total number of parameters is 570 for the whole network. After 100000 epochs of training, the $L_2$ error of the NN is reduced to the level of $\text{O}(10^{-4})$, which is not bad for discontinuous function. The result of NN based approximation of residual is shown in Fig.\ref{fig:2D_discont_NN}. We can see that the result of NN in Fig.\ref{fig:2D_discont_NN} agrees well with the residual in Fig. \ref{fig:2D_discont_rsd}. The discontinuity is well reproduced. The final results in Fig.\ref{fig:2D_discont_multiscale} also agrees well with the target function shown in Fig.\ref{fig:2D_discont_target_surf}.

\subsection{Problems of PDEs}

\subsubsection{Laplacian equation in 1D}

Consider an one-dimensional equation in nonhomogeneous form (also see \cite{Liu1995}) as follows
\begin{equation}\label{eq:1D_Laplacian}
\frac{\partial^2 u}{\partial x^2} = f(x)~,~~x \in [0,6]
\end{equation}
where
\begin{equation}
f(x) = -2 s^2 \frac{\tanh[s(x-3)]}{\sinh^2[s(x-3)]}
\end{equation}
with the essential boundary conditions
\begin{equation}
\begin{cases}
u(0) = -\tanh(3s) \\
u(6) = \tanh(3s)
\end{cases}
\end{equation}
The parameter $s$ controls the degree of localization. The analytical solution of the problem is
\begin{equation}
u(x) = \tanh (s(x-3))
\end{equation}

\begin{figure}[htbp]
	\centering
	\subfigure[Solution]{
		\label{fig:1D_Lap_u} 
		\includegraphics[width=0.7\textwidth]{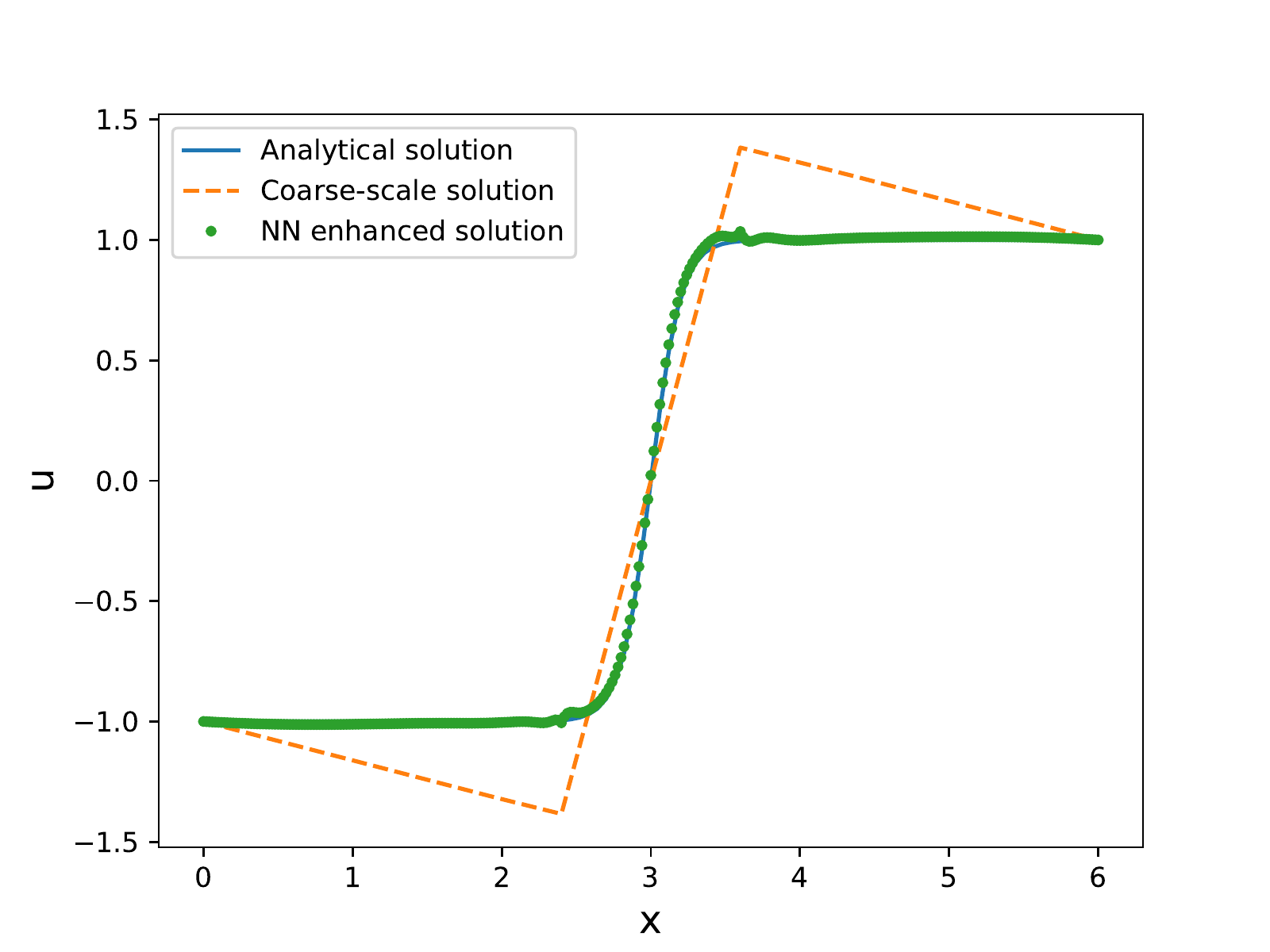}}
	\subfigure[Derivative]{
		\label{fig:1D_Lap_dudx} 
		\includegraphics[width=0.7\textwidth]{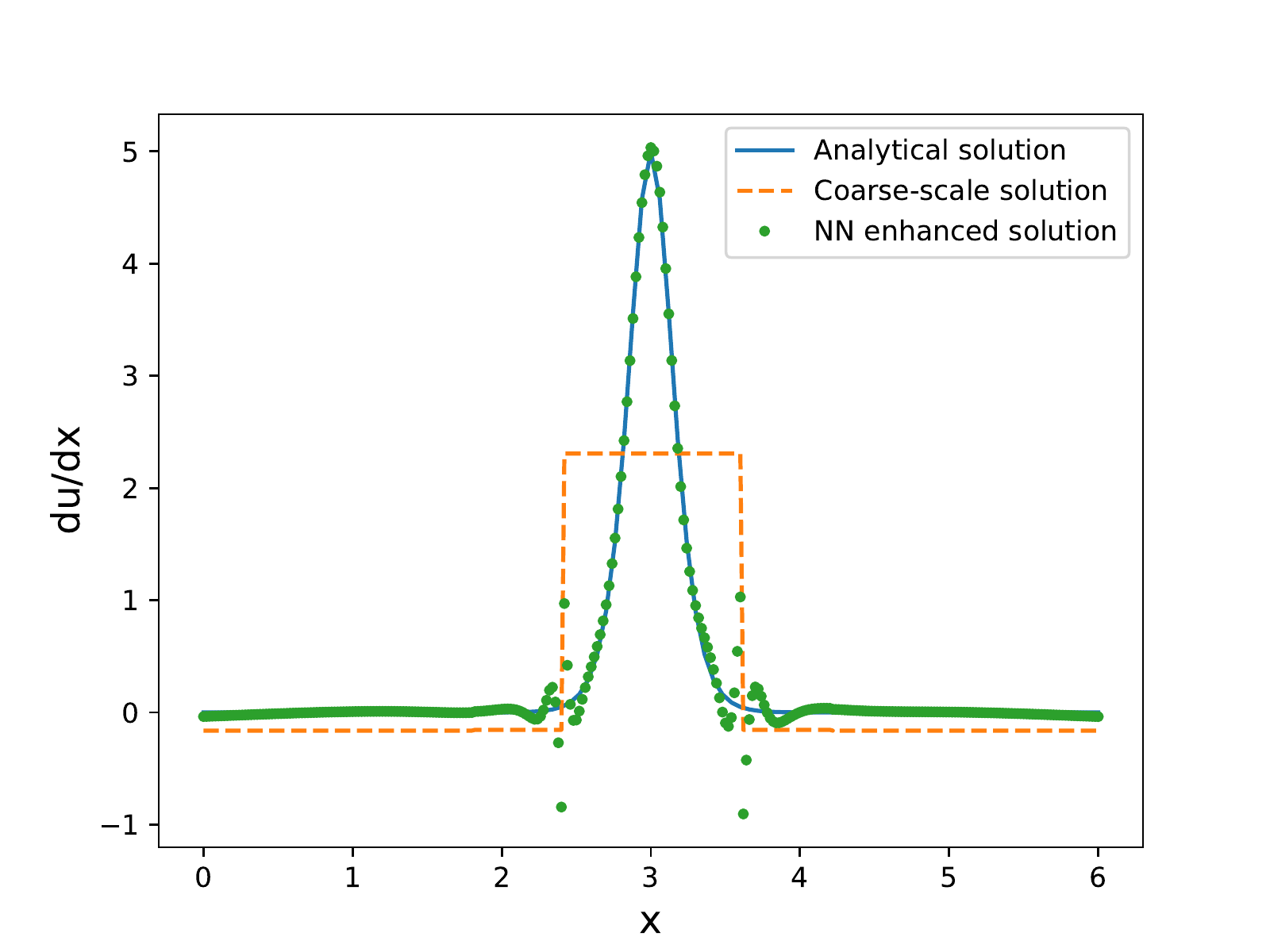}}
	\caption{\label{fig:1D_Lap} Solution of 1D Laplacian Equation (s=5)}
\end{figure}

The Galerkin form of 1D Laplacian equation in Eq.(\ref{eq:1D_Laplacian}) is not difficult to develop.
The coarse scale problem in the form of Eq.(\ref{eq:coarse_alg}) is solved based on regular finite element method. The domain$[0,6]$ is evenly divided into 10 elements. For each element, the 2-node Lagrange shape function is applied. As can be observed from Fig.\ref{fig:1D_Lap_u}, the coarse scale FEM solution is piece-wise linear. The accuracy is not good because the mesh is rather coarse. It is also not surprisingly to see in Fig.\ref{fig:1D_Lap_dudx} that the derivative solution in coarse scale is discontinuous and of lower accuracy.

By substituting the coarse-scale solution into the fine-scale governing equation Eq.(\ref{eq:fine_alg}), it could be solved based on the proposed neural networks. For the present example, a NN with three hidden layers is developed with the numbers of neurons $(4,8,5)$. The total number of parameters is 57. The activation functions are adopted as Sigmoid for the first and second hidden layers. The loss function in the energy form (see Eq.(\ref{eq:energy_loss})) is adopted. Totally 301 integration points are equally distributed within the domain $[0,6]$ for the computation of the loss function. Trained by the Adam trainer for 28000 epochs, we obtain the results in Fig. \ref{fig:1D_Lap}. It is observed that the total solution shown in Fig. \ref{fig:1D_Lap_u} agrees with the analytical solution. For the derivative solution in Fig.\ref{fig:1D_Lap_dudx}, we can see the trend of the solution is well captured by the total solution as well was the peak value. On the other hand, certain oscillations could be observed nearby the discontinuity of the coarse scale solution of  derivative. Thus it is suggested that the discontinuities of the derivative solution in coarse-scale may result in some instabilities in numerical procedure.

\begin{figure}[htbp]
	\centering
	\subfigure[Solution]{
		\label{fig:1D_Lap_u_smooth} 
		\includegraphics[width=0.7\textwidth]{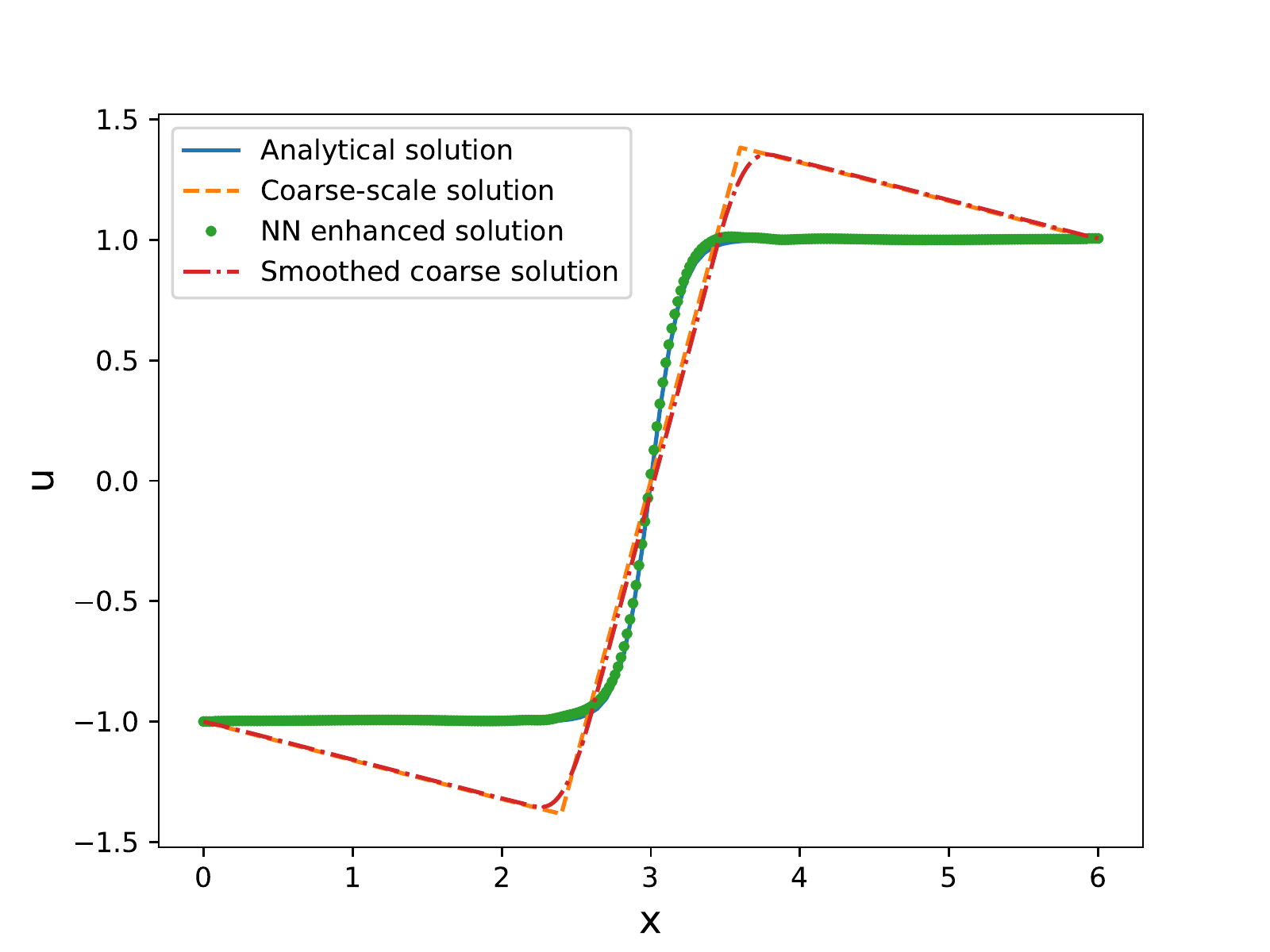}}
	\subfigure[Derivative]{
		\label{fig:1D_Lap_dudx_smooth} 
		\includegraphics[width=0.7\textwidth]{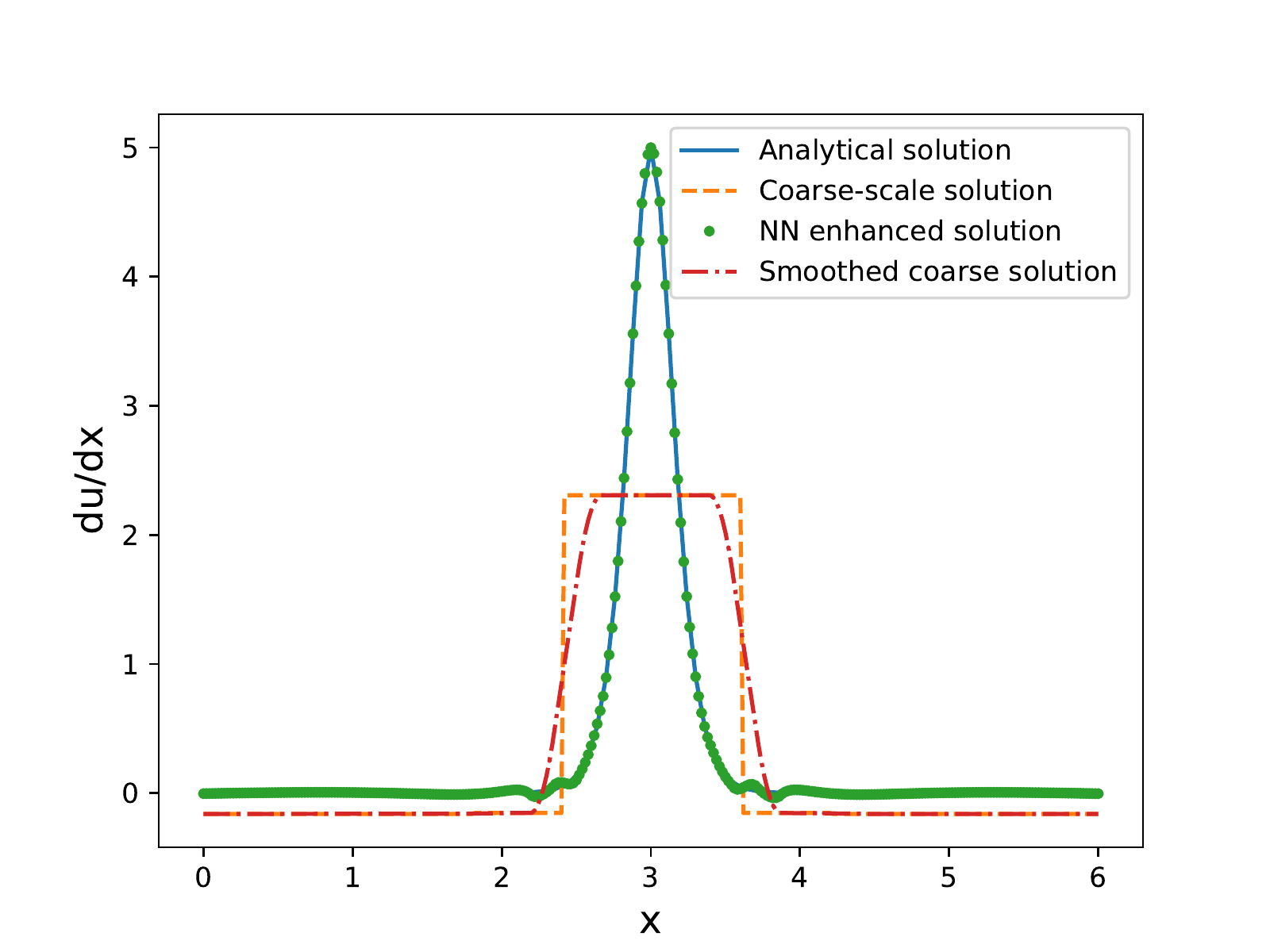}}
	\caption{\label{fig:1D_Lap_smooth} Solution of 1D Laplacian Equation with smoothing (s=5)}
\end{figure}

To overcome the numerical instability induced by the discontinuity, we propose to use the smoothed numerical solution of the coarse-scale for the numerical solution in fine-scale. As can be seen in Fig.\ref{fig:1D_Lap_smooth}, the smoothing technique well improves the solution of derivative (see Fig.\ref{fig:1D_Lap_dudx_smooth}). Moreover, the total solution in Fig.\ref{fig:1D_Lap_u_smooth} exhibits better. In the present numerical example, the smoothing method based on local averaging is used.

\begin{figure}[htbp]
	\centering
	\subfigure[Solution]{
		\label{fig:1D_Lap_u_smooth_dis} 
		\includegraphics[width=0.7\textwidth]{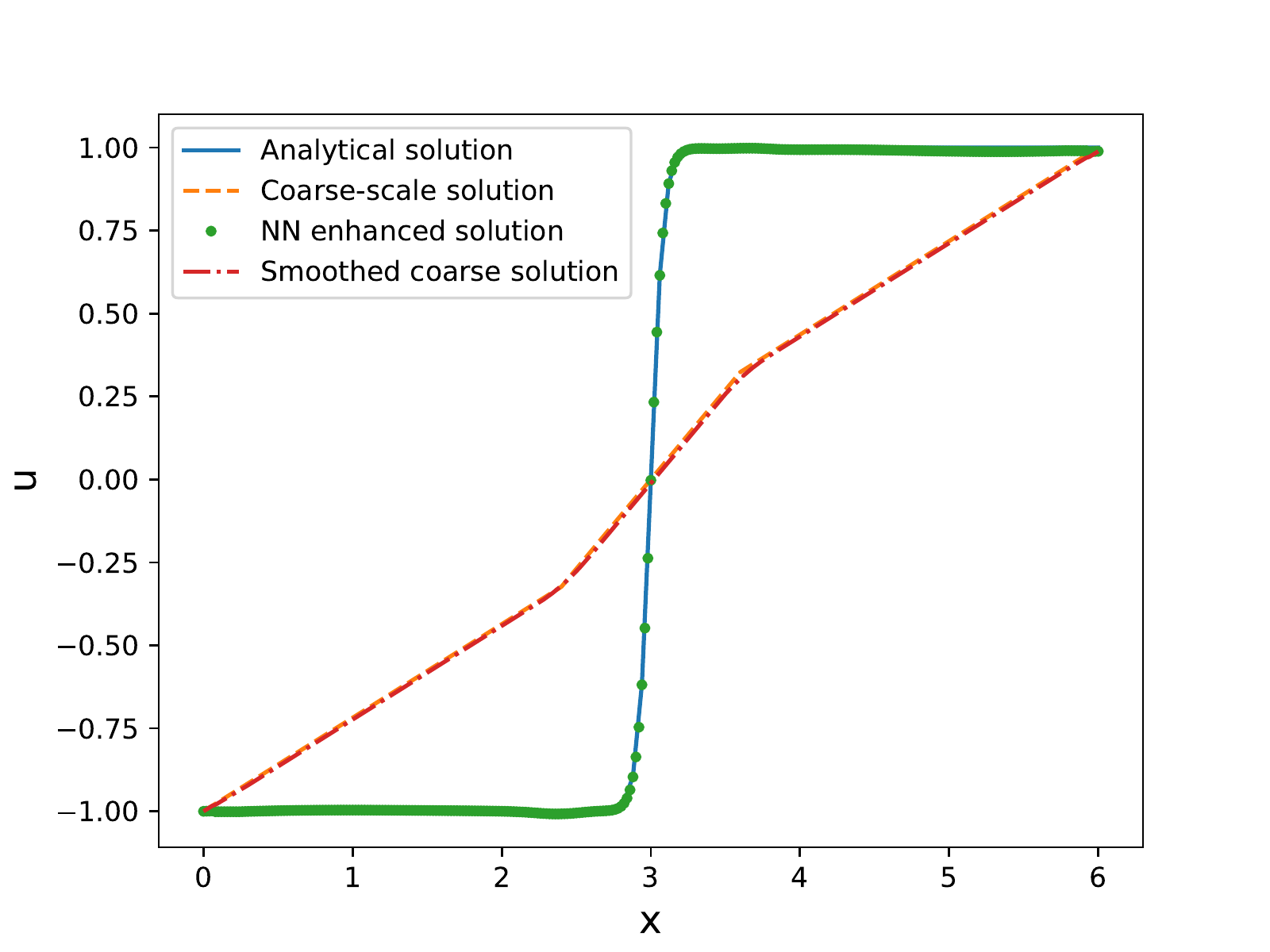}}
	\subfigure[Derivative]{
		\label{fig:1D_Lap_dudx_smooth_dis} 
		\includegraphics[width=0.7\textwidth]{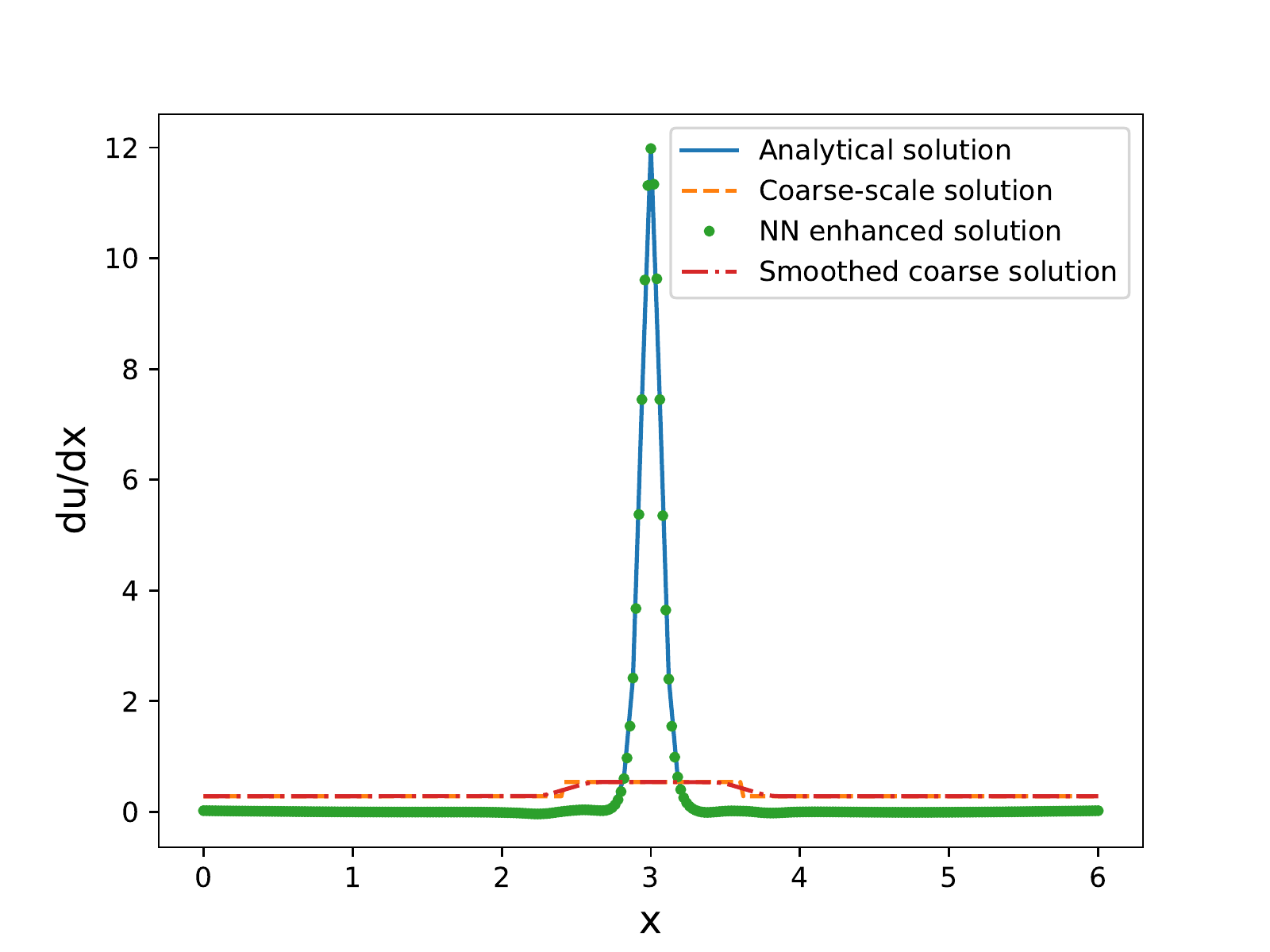}}
	\caption{\label{fig:1D_Lap_smooth_dis} Solution of 1D Laplacian Equation with smoothing (s=12)}
\end{figure}

As we know the parameter $s$ controls of the degree of localization for the derivative $\frac{\text{d} u}{\text{d} x}$. As $s$ increases, the derivative $\frac{\text{d} u}{\text{d} x}$ localizes to the center of the domain. We choose $s=5$ for the examples shown in Fig.\ref{fig:1D_Lap} and Fig.\ref{fig:1D_Lap_smooth}. Then we increase the value of $s$ from $5$ to $12$, the results are shown in Fig.\ref{fig:1D_Lap_smooth_dis}. As can be seen, the localization can hardly be reproduced by the coarse-scale solutions. On the other hand, the localization is well captured by the neural networks in fine scale. The total solutions agree well with the analytical solutions.

\subsubsection{Poisson equation in 2D}

Consider the Poisson equation in 2-D as follows
\begin{equation}
- \nabla u(\bm{x}) = f(\bm{x})~,~~\bm{x} \in \Omega
\end{equation}
where
\begin{equation}
f(\bm{x}) = 1
\end{equation}
And the essential boundary condition
\begin{equation}
u(\bm{x}) = 0 ~,~~\bm{x} \in \partial \Omega
\end{equation}
The domain is
\begin{equation}
\Omega = [-1,1]\times[-1,1] \backslash [0,1]\times[0,0]
\end{equation}

\begin{figure}[htbp]
	\centering
	\subfigure[Reference solution of $u$]{
		\label{fig:2D_Poisson_ref_u_cont} 
		\includegraphics[width=0.47\textwidth]{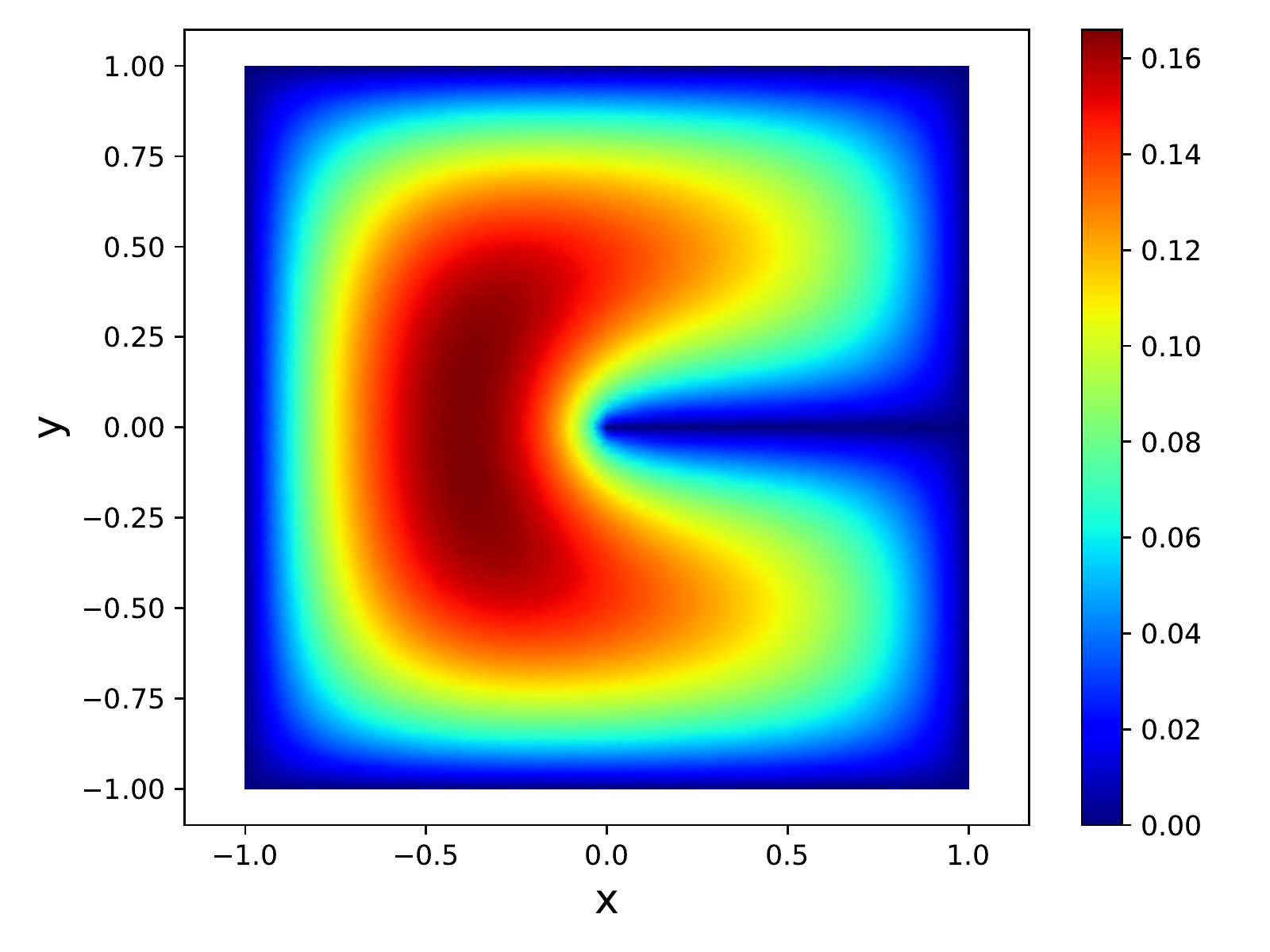}}	
	\subfigure[Reference solution of $\frac{\partial u}{\partial x}$]{
		\label{fig:2D_Poisson_ref_ux_cont} 
		\includegraphics[width=0.47\textwidth]{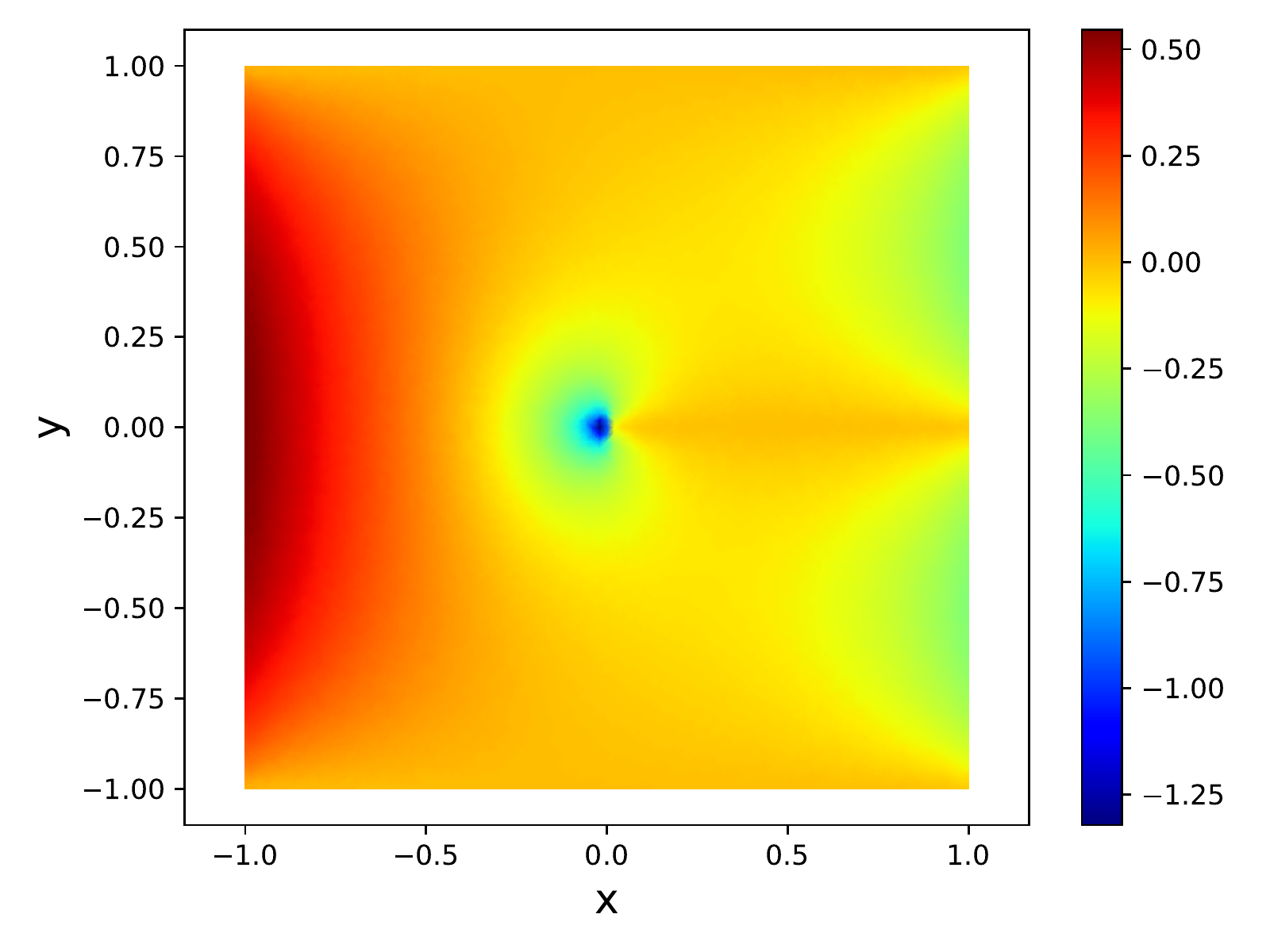}}	
	\subfigure[Reference solution of $\frac{\partial u}{\partial y}$]{
		\label{fig:2D_Poisson_ref_uy_cont} 
		\includegraphics[width=0.47\textwidth]{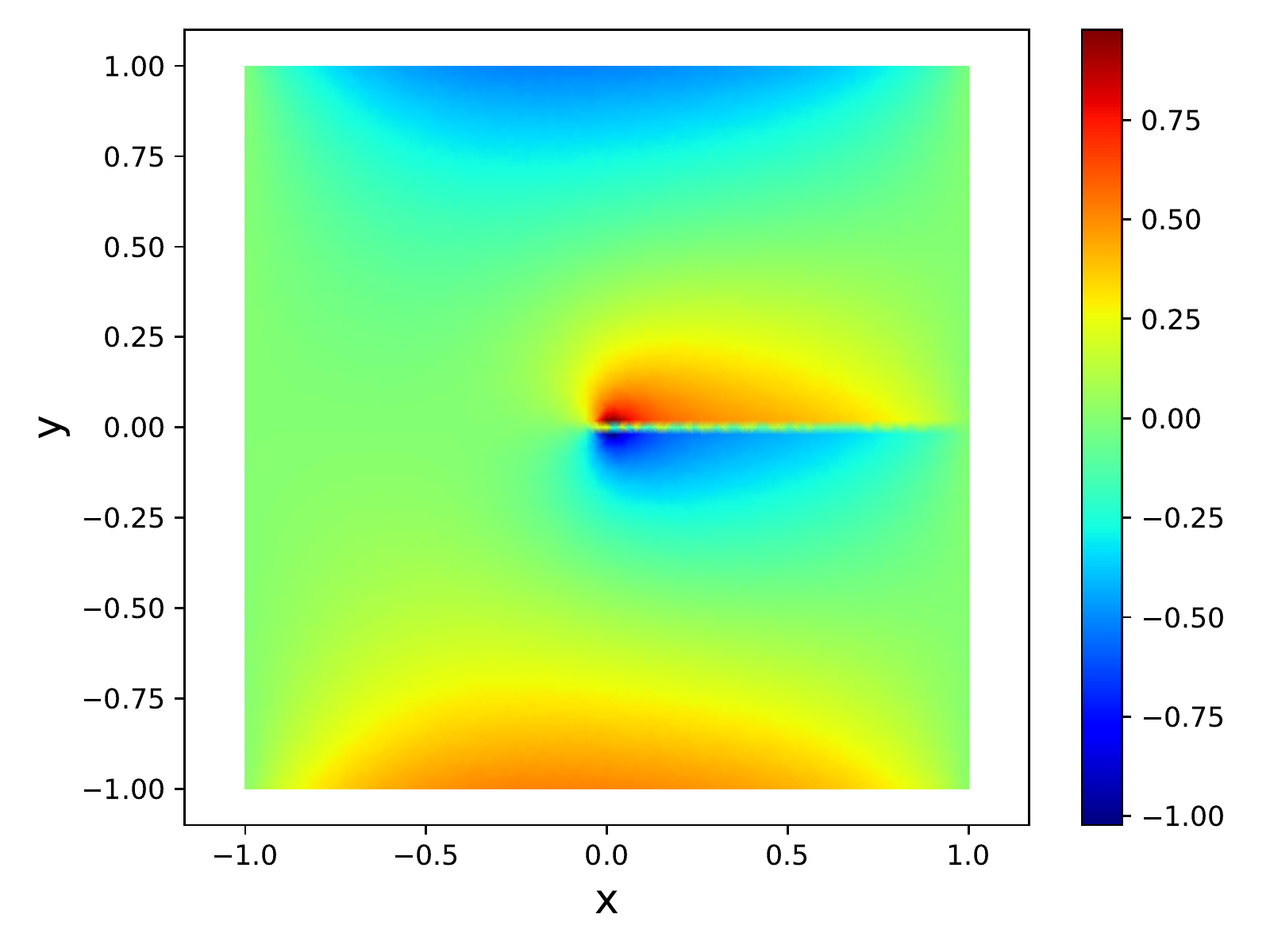}}	
	\caption{\label{fig:2D_Poisson_ref} Reference solutions based on a refined mesh (100$\times$100)}
\end{figure}

\begin{figure}[htbp]
	\centering
	\subfigure[Coarse-scale solution of $u$]{
		\label{fig:2D_Poisson_coarse_u_cont} 
		\includegraphics[width=0.47\textwidth]{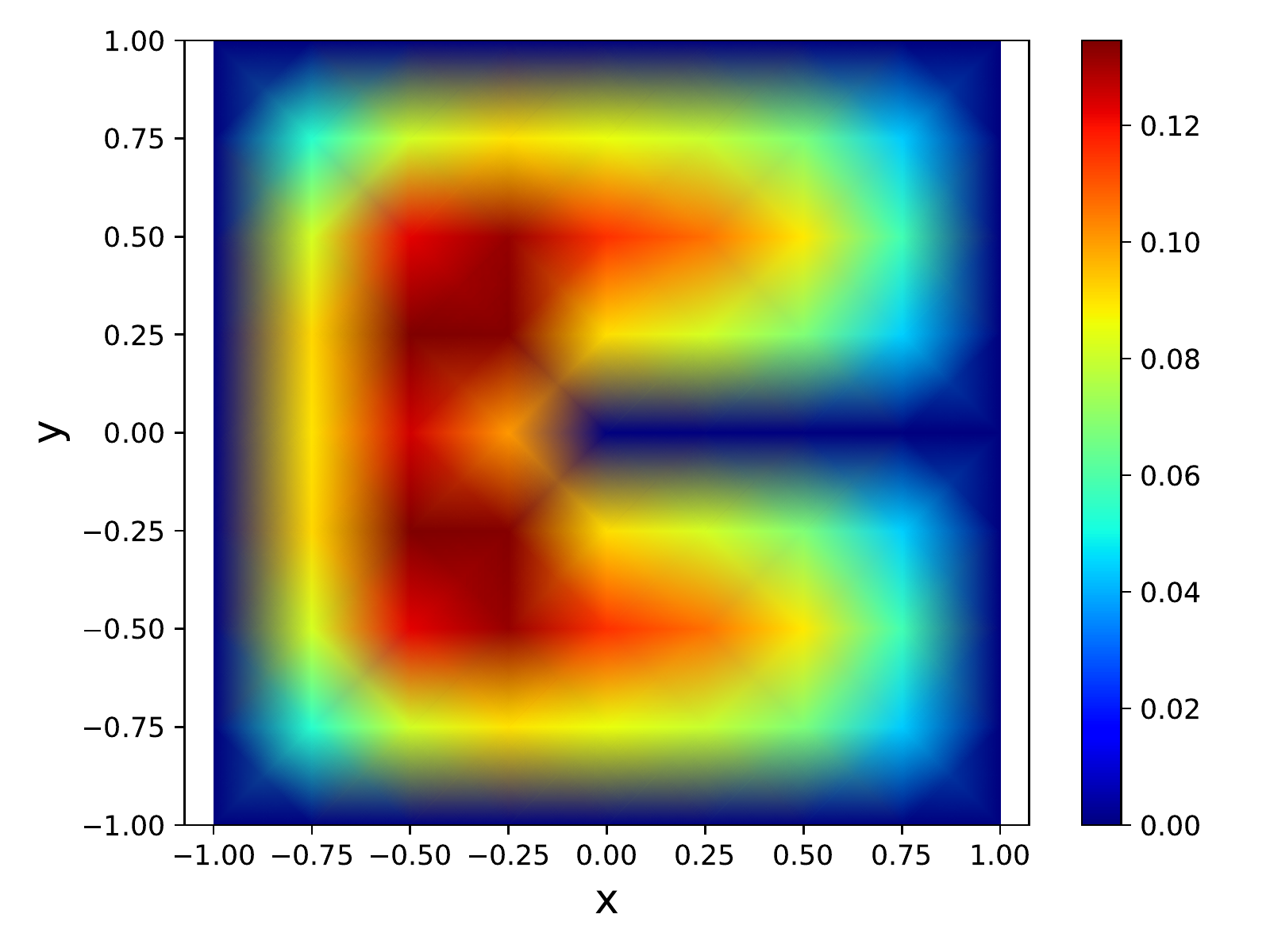}}	
	\subfigure[Coarse-scale solution of $\frac{\partial u}{\partial x}$]{
		\label{fig:2D_Poisson_coarse_ux_cont} 
		\includegraphics[width=0.47\textwidth]{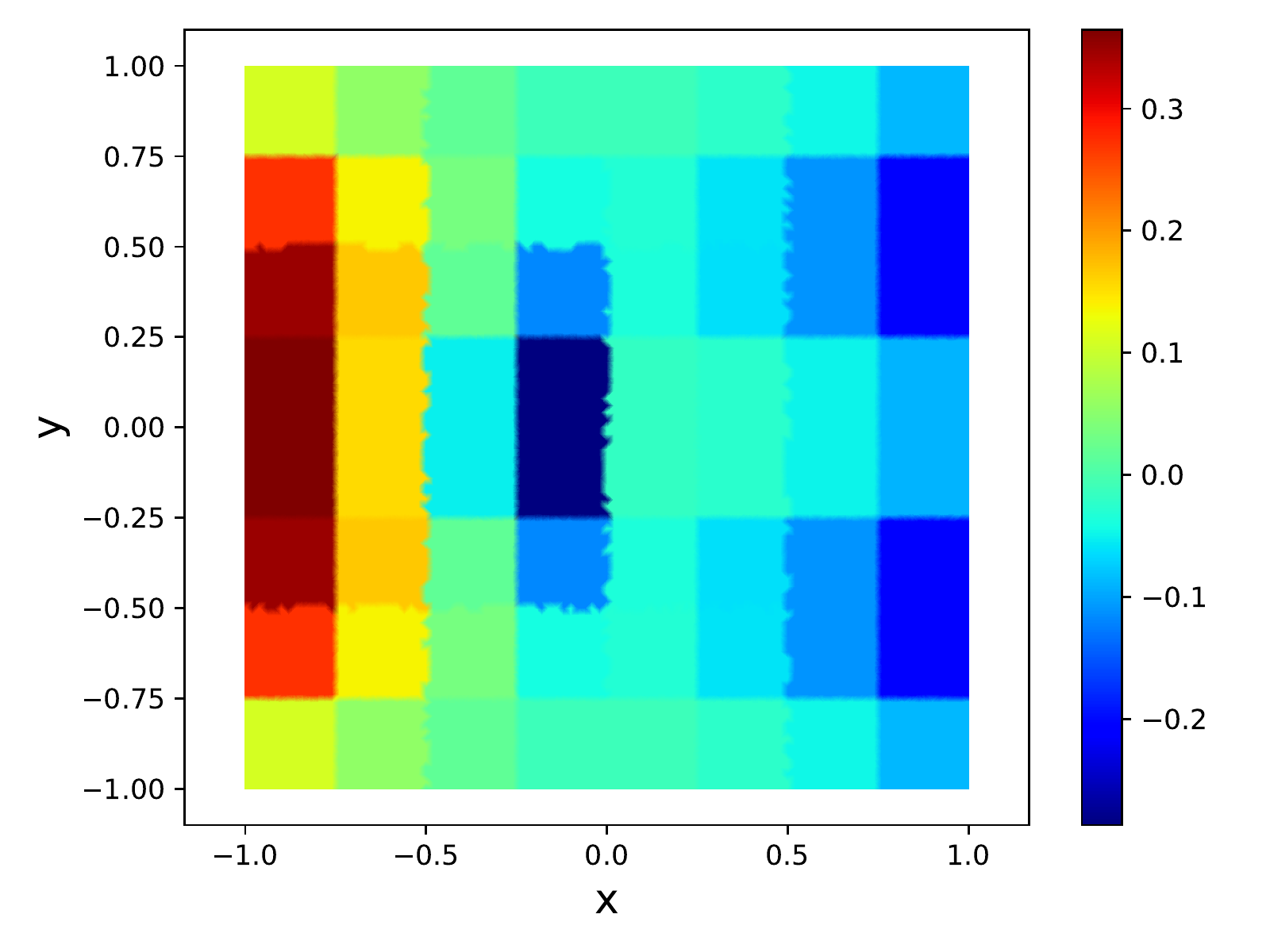}}	
	\subfigure[Coarse-scale solution of $\frac{\partial u}{\partial y}$]{
		\label{fig:2D_Poisson_coarse_uy_cont} 
		\includegraphics[width=0.47\textwidth]{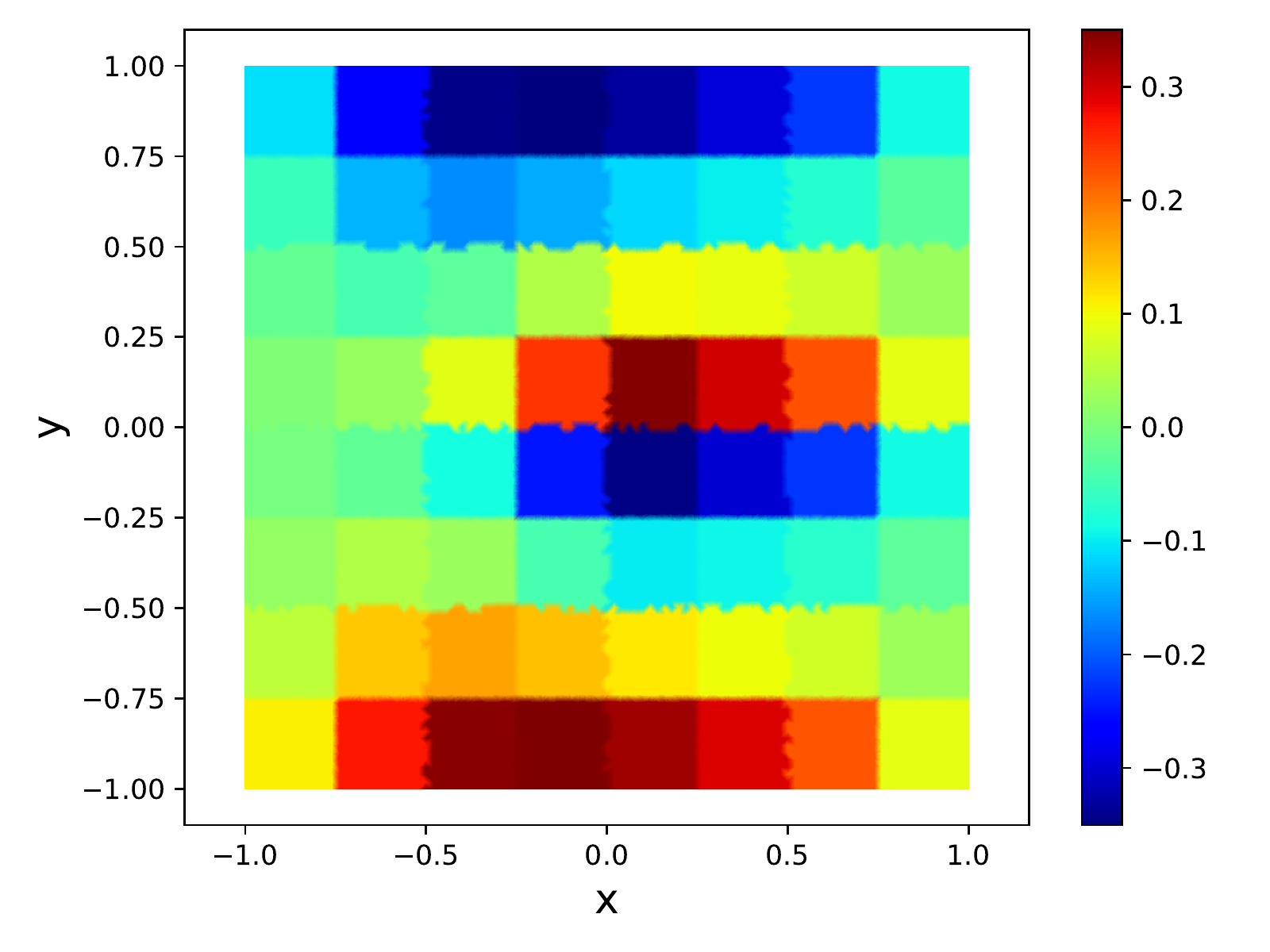}}	
	\caption{\label{fig:2D_Poisson_corase} Coarse-scale solutions based on a coarse mesh (8$\times$8)}
\end{figure}

This problem is widely used to illustrate the methods with self-adaption or multi-scale nature. By using a rather refined mesh (100$\times$100), the reference solution (also see \cite{EWeinan2018}) of this problem is shown in Figure \ref{fig:2D_Poisson_ref}. The bi-linear rectangular element with reduced integration and hourglass control \cite{Belytschko2013} is used for the simulation. The solution of $u$ shown in Fig. \ref{fig:2D_Poisson_ref_u_cont} is continuous and smooth, without any hourglass patterns. The singularity at center (see Fig. \ref{fig:2D_Poisson_ref_ux_cont}) could be clearly observed for the solution of $\frac{\partial u}{\partial x}$. In the coarse-scale, a rather coarse mesh, i.e. 8$\times$8, is developed. The results of the coarse-scale solution are shown in Figure \ref{fig:2D_Poisson_corase}. It is not surprisingly to see that the resolution of the coarse-scale solution is relatively low. Also the patch-wise discontinuity in the fields of gradients (i.e. $\frac{\partial u}{\partial x}$ and $\frac{\partial u}{\partial y}$) could be clearly observed. The reason is that the finite element base functions used for the coarse-scale solution is of $C_0$ continuity.

\begin{figure}[htbp]
	\centering
	\subfigure[NN solution of $\tilde{u}$ in fine-scale]{
		\label{fig:2D_Poisson_fine_u_cont} 
		\includegraphics[width=0.43\textwidth]{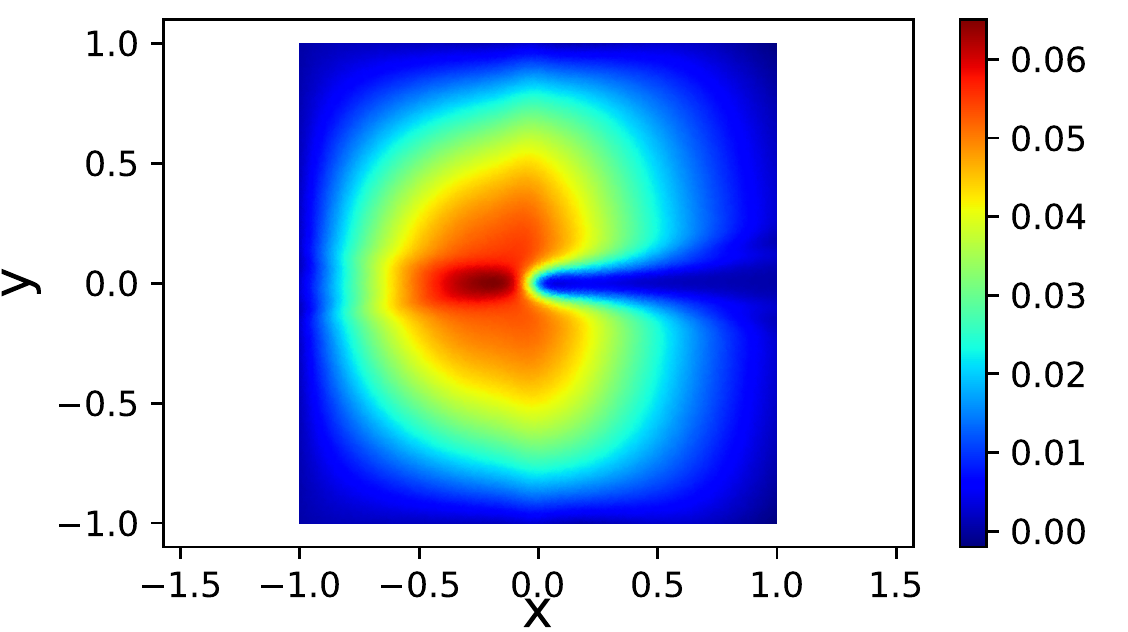}}	
	\subfigure[Multi-scale solution of $u$]{
		\label{fig:2D_Poisson_total_u_cont} 
		\includegraphics[width=0.4\textwidth]{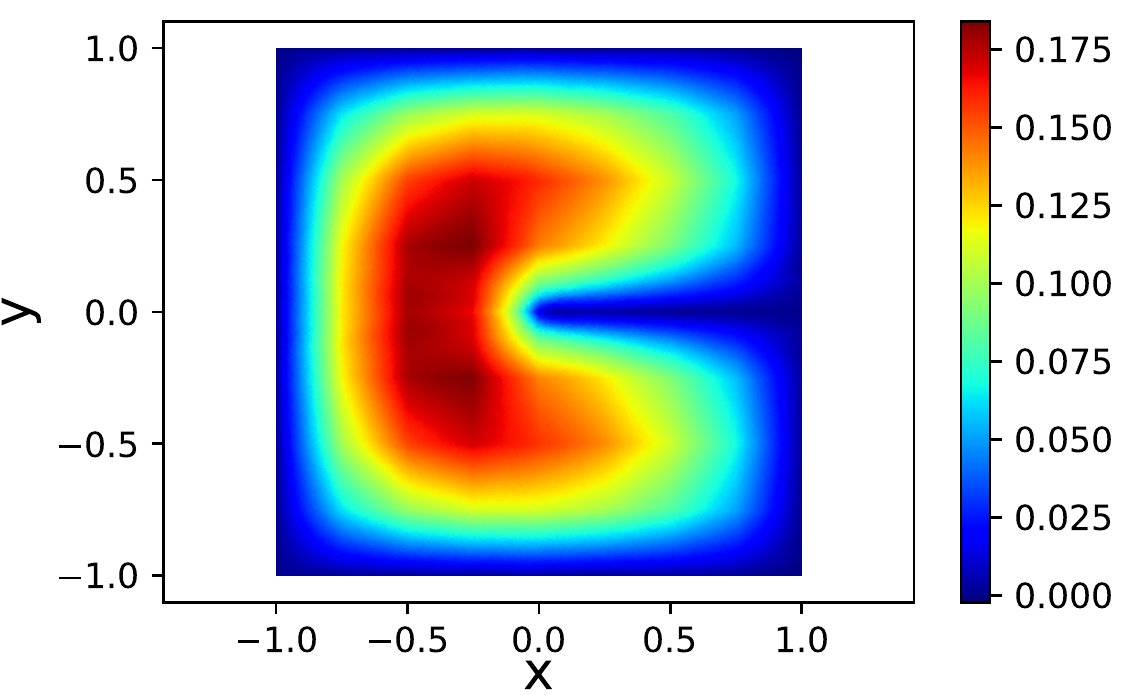}}	
	\subfigure[Coarse-scale solution of $\frac{\partial u}{\partial x}$]{
		\label{fig:2D_Poisson_total_ux_cont} 
		\includegraphics[width=0.4\textwidth]{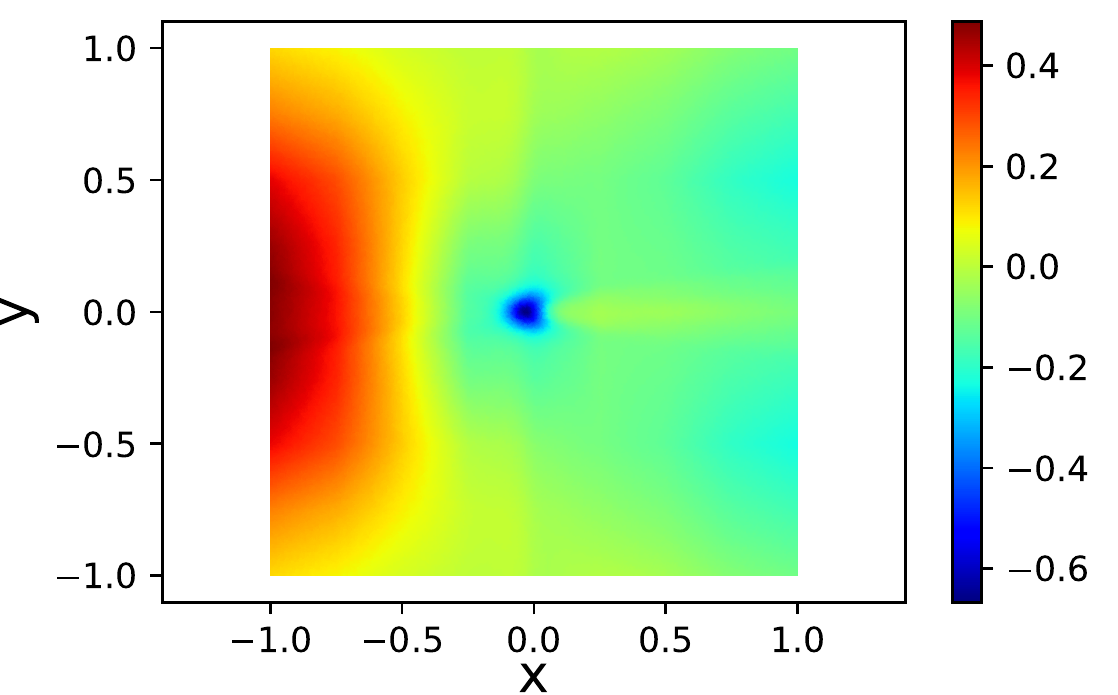}}	
	\subfigure[Multi-scale solution of $\frac{\partial u}{\partial y}$]{
		\label{fig:2D_Poisson_total_uy_cont} 
		\includegraphics[width=0.4\textwidth]{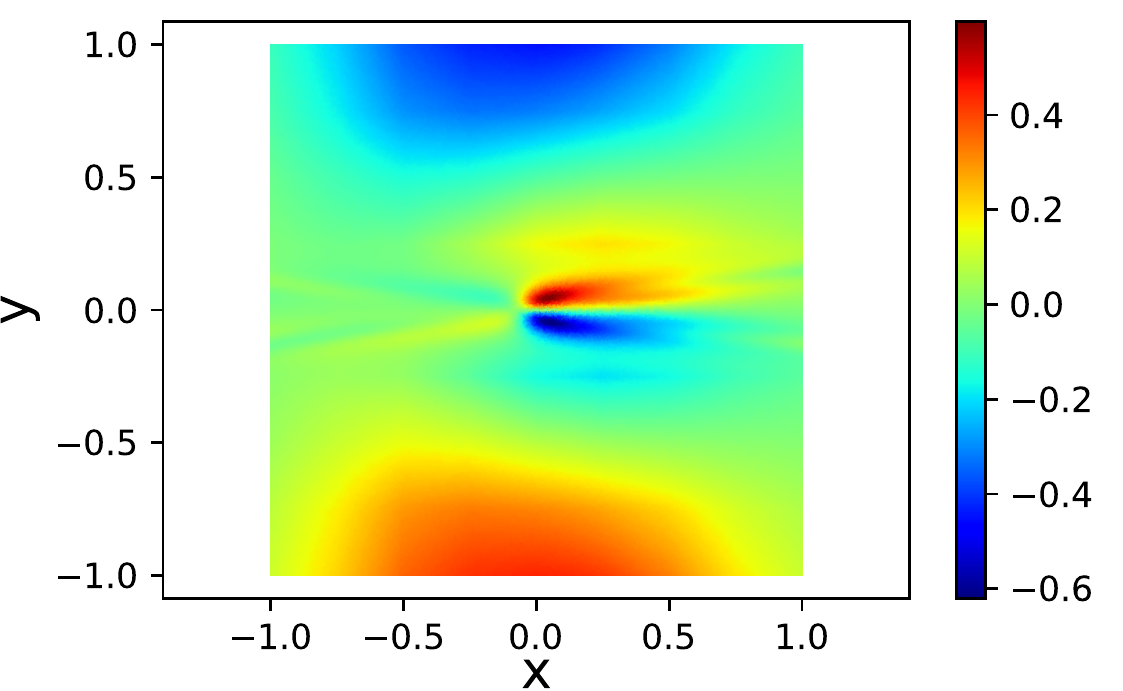}}	
	\caption{\label{fig:2D_Poisson_total} Multi-scale solutions with NN enhancement(Case 1)}
\end{figure}

\begin{figure}[htbp]
	\centering
	\subfigure[NN solution of $\tilde{u}$ in fine-scale]{
		\label{fig:2D_Poisson_fine_u_cont_1} 
		\includegraphics[width=0.42\textwidth]{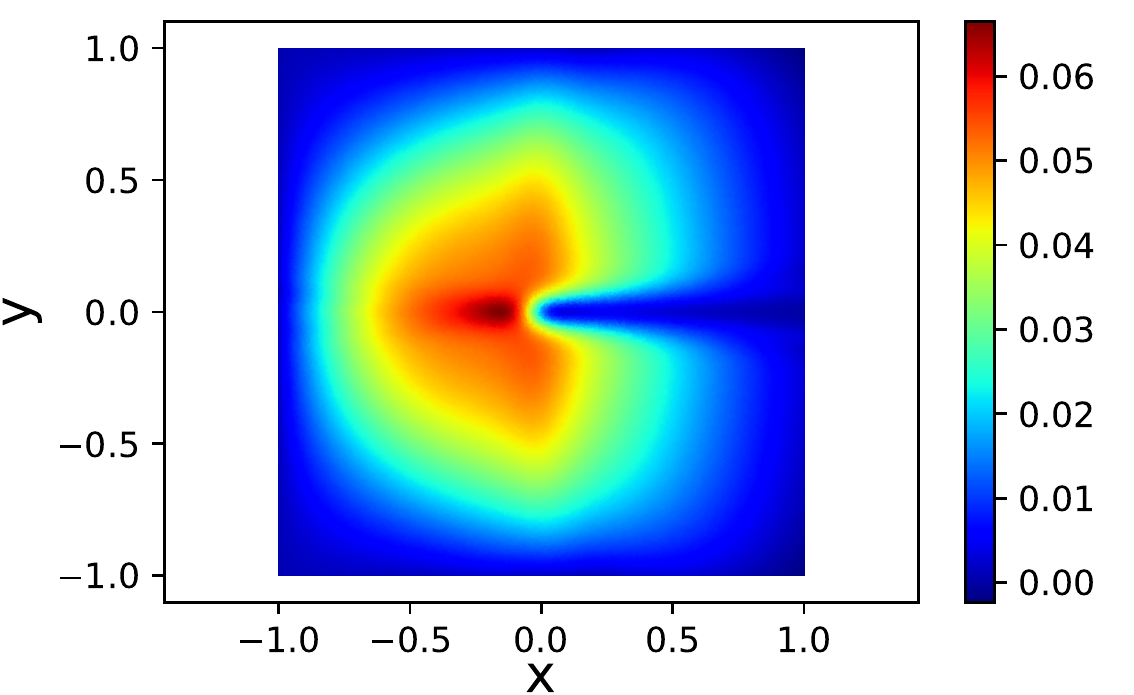}}	
	\subfigure[Multi-scale solution of $u$]{
		\label{fig:2D_Poisson_total_u_cont_1} 
		\includegraphics[width=0.42\textwidth]{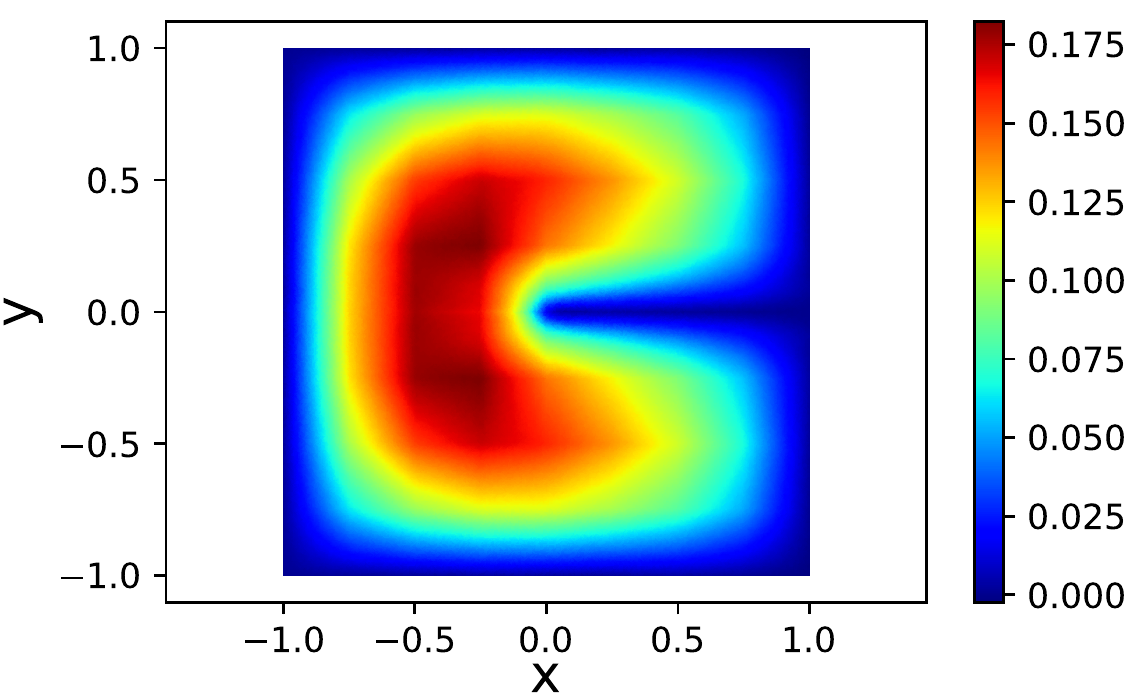}}	
	\subfigure[Coarse-scale solution of $\frac{\partial u}{\partial x}$]{
		\label{fig:2D_Poisson_total_ux_cont_1} 
		\includegraphics[width=0.42\textwidth]{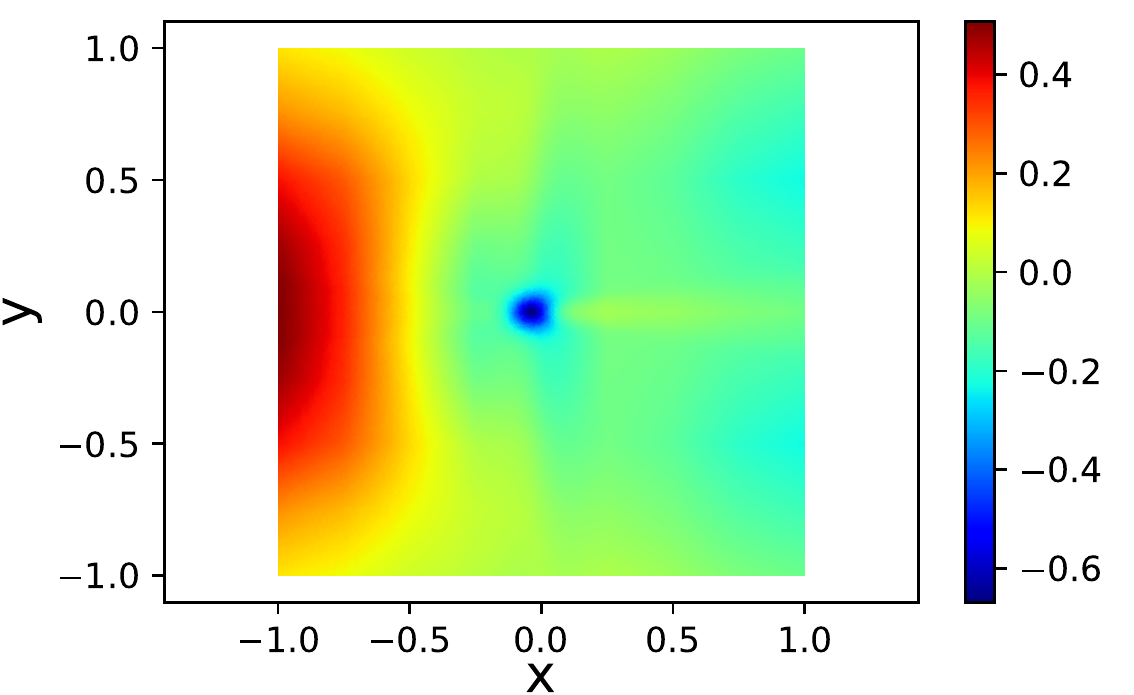}}	
	\subfigure[Multi-scale solution of $\frac{\partial u}{\partial y}$]{
		\label{fig:2D_Poisson_total_uy_cont_1} 
		\includegraphics[width=0.39\textwidth]{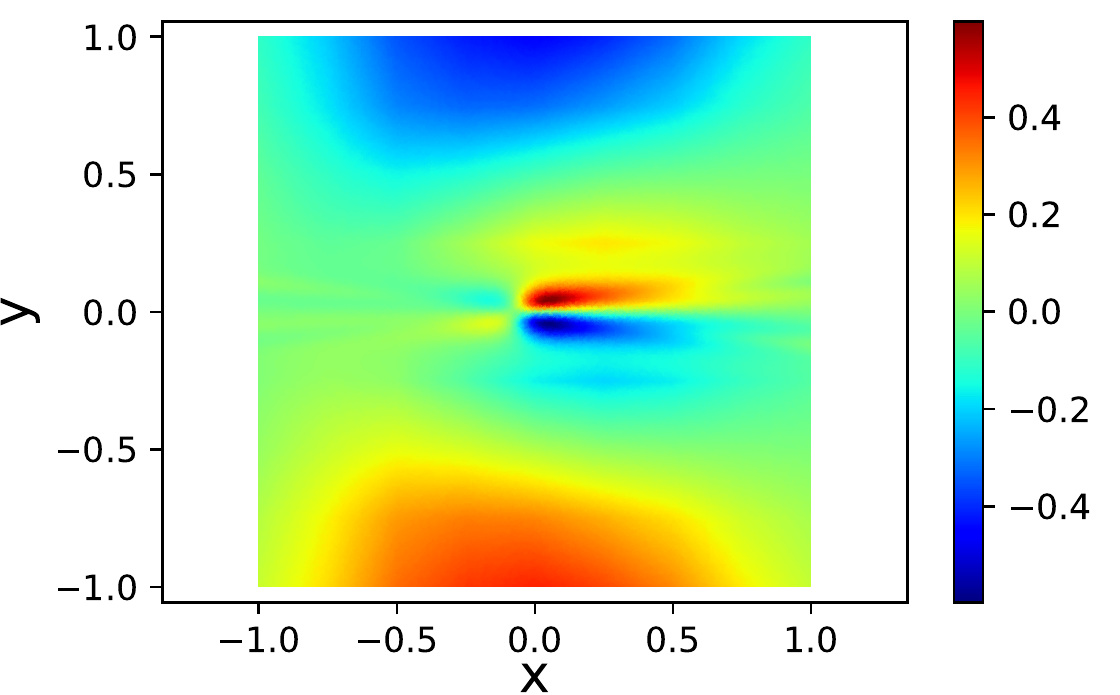}}	
	\caption{\label{fig:2D_Poisson_total_1} Multi-scale solutions with NN enhancement (Case 2)}
\end{figure}

Based on the proposed multi-scale method, the governing equation could be solved based on NN with the custom loss function in the form of Eq.(\ref{eq:energy_loss}), for which the influence from the coarse-scale solution is also taken into account. To suppress the possible instability in the computation at fine scale, a regular smoothing is performed to the coarse-scale solution. A neural network of 5 densely connected hidden layers is developed. The numbers of neurons in the layers are $(10,15,25,15,10)$. The Sigmoid activation function is chosen for the neurons in hidden layers. The total number of unknowns for the NN is 1276. As we can see, the neural network we used here is not special in anyway. To calculate the loss function in the integration form, totally $151\times151$ integration points are assigned to the domain $[-1,1]\times[-1,1]$. The essential boundaries are also covered by the set of integration points. The NN was trained for numbers of times, the results of two cases are shown in Fig.\ref{fig:2D_Poisson_total} and Fig.\ref{fig:2D_Poisson_total_1}. 
Each cases is trained by Adam trainer for 220,000 epochs. The multi-scale solutions of both cases shown in Figures \ref{fig:2D_Poisson_total} and \ref{fig:2D_Poisson_total_1} agree with the reference solution in Fig.\ref{fig:2D_Poisson_ref}. The accuracy of the multi-scale solutions is well improved comparing to the coarse-scale solution in Fig.\ref{fig:2D_Poisson_corase}. Especially, Figures \ref{fig:2D_Poisson_ref_ux_cont} and \ref{fig:2D_Poisson_ref_ux_cont} indicates that the singularity for the gradient field (i.e. $\frac{\partial u}{\partial x}$) is captured by the multi-scale solution. On the other side, the results of NN enhancements for both cases shown in Fig.\ref{fig:2D_Poisson_fine_u_cont} and Fig.\ref{fig:2D_Poisson_fine_u_cont_1} are not very close to each other, although they are kind of similar. These results indicate that some stable training techniques, for example BFGS methods \cite{Matthies1979,Fletcher2000,Ren2018}, are needed and deserve more efforts in the future work. 

\section{Concluding remarks}
\label{sec:VI}
The classic methods for PDEs (e.g. finite element, finite difference, meshfree methods and so on) have been developed in more than 50 years. They are mature and well-equipped with robustness and efficiency. The NN/AI based methods proposed in recent years are rather different, but not the replacements of them ( see also in\cite{RAISSI2019686}). In this context, one may ask: are there any approaches to make good use both of them? The present work could be considered as a trial.

By introducing the multi-scale concept, the classic methods and the NN based methods are adopted to solve the problems in coarse-scale and fine-scale, respectively. Then they are capable of working within a unified framework in harmony. The implementation of the proposed methods are also not difficult because the interactions between scales are weak. That is to say, no iteration between scales is required in the present work. The proposed methods offer the possibility to enhance the capacity of the classic methods by using NN based methods, which are at present undergoing rapid developments made by world-class communities.

In addition, there are still interesting and important issues. How to make the choices of the structure, activation function and training strategy of the NN with consideration of influence from the other scale? How to analyze the function space with linear and nonlinear constructs? Are there any preferable or optimal combinations of classic methods and NN based methods? The present work is just opening a door to the subject.

%%\section*{Acknowledgments}
%%This was was supported in part by......

\appendix
\section*{Appendix A: Schemes for numerical integration}
\label{A:1}
For the functional norms adopted in the loss function, methods of numerical integration are usually introduced in the implementation. Sometimes, the methods of numerical integration result in rather different accuracy and stability for the whole procedure of numerical simulation. 

\begin{itemize}
	\item Cell integration \\
	The integration domain is divided into a finite number of cells.
	\begin{equation}
	\mathcal{I}(g) = 
	\int_{\Omega} g(\bm{x}) \text{d} \Omega
	= \sum_{I = 1}^{N_\text{C}} \int_{\Omega_I} g(\bm{x}) \text{d} \Omega
	\end{equation}
	The integration over the $I$-th cell $\Omega_I$ could be approximated by certain quadrature rules, for example Gauss quadrature rule, as follows
	\begin{equation}
	\mathcal{I}_I(g) =
	\int_{\Omega_I} g(\bm{x}) \text{d} \Omega 
	\approx
	\mathcal{I}^\text{h}_I(g) =
	\sum_{J = 1}^{N_{\text{Int}}} \omega_J g (\bm{\xi}_J) \det{\mathcal{J}}
	\end{equation}
	where $\omega_J$ is the $J$-th Gaussian weight; and the Jacobian matrix is
	\begin{equation}
	\mathcal{J} = \frac{\partial \bm{x}}{\partial \bm{\xi}}
	\end{equation}
	The precision of the cell based integration is usually not bad if the function $g(\cdot)$ is of good smoothness. The error could be easily estimated from the error of integration in each cell based on classic theory of quadrature. As we can see, the cell-based quadrature rule works well with the mesh-based methods, say finite element methods.
	
	\item Nodal integration \\
	For the particle-based methods without background mesh, the methods of nodal integration \cite{BEISSEL199649,Chen2001} are usually adopted. We can see the following form:
	\begin{equation}
	\mathcal{I}(g) =
	\int_{\Omega} g (\bm{x}) \text{d} \Omega
	\approx 
	\mathcal{I}^\text{h}(g) =
	\sum_{I = 1}^{N_\text{P}} g (\bm{x}_I) V_I
	\end{equation}
	where $V_I$ is the representative volume of the node $\bm{x}_I$. The selection of the nodal set $\{\bm{x}_1,...,\bm{x}_2,...,\bm{x}_I,...,\bm{x}_{N_\text{P}} \}$ is a elaborate work, especially for high dimensional problems. For the grid-based nodal set, the error of the integration reads
	\begin{equation}
	\mathcal{I}(g) - \mathcal{I}^\text{h}(g) \simeq
	\frac{\mathcal{C}(g)}{(N_\text{P})^{\alpha / d}}
	\end{equation}	
	where $\mathcal{C}(g)$ is a constant related to function $g(\cdot)$; $\alpha$ is related to the smoothness of $g(\cdot)$; and $d$ is the dimension of the problem. It could be easily observed that the error is strongly dependent to the dimension of the problem, which is demoted by $d$.
	
	In recent years, the random or pseudo-random data sets have drawn much attention from the community. By selecting the nodal set $\{\bm{x}_1,...,\bm{x}_2,...,\bm{x}_I,...,\bm{x}_{N_\text{P}} \}$ via random or pseudo-random sampling, the error norm could be established based on expectation operator $\mathbf{E}$ as follows:
	\begin{equation}\label{eq:MC_I}
	\mathbf{E} [\mathcal{I}(g) - \mathcal{I}^\text{h}(g)]^2
	= \frac{\mathbf{var}(g)}{N_\text{P}}
	\end{equation}
	where the variance operator is
	\begin{equation}
	\mathbf{var}(g) = \mathbf{E}g^2 - (\mathbf{E}g)^2
	\end{equation}
	The dimensionless convergence rate could be achieved according to Eq.(\ref{eq:MC_I}) (see \cite{SLOAN19981,dick_kuo_sloan_2013}). Nevertheless, the convergence is random. On the other side, the quadrature based on random or pseudo-random data set works well in the framework of neural networks\cite{Han2018,hennig2022}.	
	
\end{itemize}

%Bibliography
\bibliographystyle{unsrt}  
\bibliography{ref}

\end{document}